\documentclass[12pt]{article}
\usepackage[english]{babel}
\usepackage[utf8]{inputenc}
\usepackage{amsthm,amsmath,amsfonts,amssymb}
\usepackage{bbold}
\usepackage[normalem]{ulem}
\usepackage{nicefrac}
\usepackage{url}
\usepackage{mathtools}
\usepackage{color}
\usepackage{anysize}
\usepackage{geometry}
\usepackage{comment}
\usepackage{dsfont}
\usepackage{enumitem}
\marginsize{2cm}{2cm}{1,5cm}{1,5cm}

\numberwithin{equation}{section}
\theoremstyle{plain}

\newtheorem{theorem}{Theorem}[section]
\newtheorem{lemma}[theorem]{Lemma}
\newtheorem{proposition}[theorem]{Proposition}
\newtheorem{remark}[theorem]{Remark}
\newtheorem{corollary}[theorem]{Corollary}
\newtheorem{definition}[theorem]{Definition}

\usepackage{amsmath}
\usepackage{array}
\usepackage{graphicx}

\graphicspath{{./Images/}}
\usepackage{booktabs}
\newcommand{\1}{\mathds{1}}

\mathtoolsset{showonlyrefs}
\usepackage[backend=bibtex,style=numeric]{biblatex}
\definecolor{myblue}{RGB}{20,60,160}
\definecolor{myred}{RGB}{200,16,40}
\usepackage[colorlinks = true,
linkcolor = myblue,
citecolor = myred 
]{hyperref}

\title{Probabilistic counting of labelled trees with prescribed vertical edge profile}

\author{Vanessa Dan and Alexis Metz-Donnadieu}

\begin{document}
\maketitle
\begin{abstract}
The vertical edge profile of a labelled plane tree $t$ is the collection of integers counting the number of oriented edges from label $i$ to label $j$ for every possible pairs of labels $(i, j)$. In \cite{MarkovPropertyMD}, it was shown that a process constructed from the vertical edge profile forms a Markov chain for a wide class of models of random labelled plane trees where edge label increments are constrained to lie in $\{-1, 0, 1\}$. The purpose of this paper is twofold. First, we continue the study of Markov properties of the vertical edge profile. In particular, we derive general formulas for the transition kernel and the initial distribution of the class of chains considered in \cite{MarkovPropertyMD}. Second, we apply these results to the probabilistic enumeration of trees with a prescribed vertical edge profile. These applications yield a unified probabilistic derivation of several enumeration formulas. In particular, we recover some results established in \cite{BousquetChapuy} for binary trees and Cayley trees.
These enumeration formulas admit a remarkable factorised form, already noted in \cite{Bernardi, BousquetChapuy}. The probabilistic methods developed here explain this phenomenon by showing that these factorizations arise naturally from the Markov structure established in \cite{MarkovPropertyMD}.
\end{abstract}
\section{Introduction}
 We study random plane trees equipped with an integer-valued vertex labelling whose edge increments belong to $\{-1, 0, 1\}$. In probability theory, these objects arise in a variety of contexts, from the study of branching processes and particle systems \cite{Chauvin,DevroyeGW, Janson}, to their connections with planar maps and models of random geometry \cite{CoriVauquelin,BrownianMapMiermont, Schaeffer}. They have also been extensively studied in combinatorics \cite{Bernardi, BousquetMelou,BousquetChapuy}. Establishing connections between combinatorial and probabilistic methods has long been a central theme in the study of random trees. Such connections arise at many levels, ranging from enumeration techniques to methods in asymptotic probability, and provide one of the main motivations for the present work.
 
 \medskip
 
 The class of models we consider encompasses branching random walks as well as combinatorial models such as random embedded binary and ternary trees, and Cayley or plane trees with uniformly distributed edge displacements. More precisely, it is defined as follows.
  Fix a non-trivial critical or subcritical distribution $\xi$ on $\{0, 1, 2, \cdots \}$ and a family $(\eta^{(d)}; d\geq 1)$, where for every $d\geq 1$, $\eta^{(d)}$ is a probability distribution on $\{-1, 0, 1\}^d$. We first sample a Bienaymé--Galton--Watson tree $T$ with offspring distribution $\xi$, viewed as a random plane tree. Then, conditionally on $T$, we sample independently for every non-leaf vertex $v$ of $T$ a random vector $L_v\in \{-1, 0, 1\}^{c_T(v)}$ with law $\eta^{(c_T(v))}$, where $c_T(v)$ denotes the number of children of $v$ in $T$. We then label the edges connecting $v$ to its children using the $c_T(v)$ components of $L_v$ according to the left-to-right order induced by the plane embedding. The label $\ell_T(u)$ of a vertex $u$ of $T$ is defined as the sum of the edge labels along the ancestral path from the root (labelled $0$ by convention) to $u$.
  
  \medskip

The first objective of this paper is to continue the study of a process introduced in \cite{MarkovPropertyMD}, called the vertical edge profile. The vertical edge profile of a labelled tree $t$ is defined as the family $(e^t_{i, j})_{i,j \in \mathbb Z, |i-j|\leq 1}$, where $e^t_{i, j}$ denotes the number of edges in $t$ from a parent with label $i$ to a child with label $j$. It was shown in \cite{MarkovPropertyMD} that the process
\begin{align*}
	(e^T_{k-1, k}, e^T_{k, k-1})_{k\geq 1},
\end{align*}
is a Markov chain when the labelled tree $(T, \ell_T)$ is sampled as described above. In this work, we consider the process
\begin{align*}
		(e^T_{k-1, k}, e^T_{k-1, k-1}, e^T_{k, k-1})_{k\geq 1},
\end{align*}
which remains a Markov chain (Proposition \ref{MP}). We provide explicit formulas for the initial distribution and transition kernel of this Markov chain (Theorem \ref{FormKernel} and Theorem \ref{th_profat0}). The formulas are expressed in terms of convolutions of random variables naturally associated with the offspring and displacement distributions. This is in keeping with the classical correspondence between Bienaymé--Galton--Watson trees and random walks. In particular, encodings such as the Łukasiewicz walk play a central role throughout this work.
\medskip

The second objective of this paper is to apply these formulas to the enumeration of trees in some combinatorial classes with prescribed vertical edge profile.
Probabilistic approaches to the enumeration of trees have a long history, going back to the study of Bienaymé--Galton–Watson branching processes and their connections with combinatorial tree models. In particular, many classical families of trees can be realized as Bienaymé--Galton–Watson trees that are conditioned on their size (see e.g. \cite{Janson}), thereby reducing enumeration problems to the computation of probabilities in suitable random tree models. The key observation underlying these methods is that, for suitable choices of offspring and displacement distributions, the induced random tree is uniformly distributed over a combinatorially interesting class. In particular, if we wish to enumerate a given subclass of trees, it suffices to compute the probability that the random tree lies in this subclass. In the models considered here, the formulas for the transition kernel and the initial distribution of the Markov chain allow us to compute the probability that a tree has a given vertical edge profile (see Proposition \ref{Path_probabilities}). As an application, we give enumeration formulas (some of which were previously proved by different methods, see in particular \cite{Bernardi, BousquetChapuy}) for the number of trees with a prescribed vertical edge profile $(x_{i, j})_{|i-j|\leq 1}$ in the following classes:
\begin{itemize}[label=\scriptsize$\bullet$]
	\item Theorem \ref{theorem_embeded}.2 : rooted labelled Cayley trees with label increments along edges constrained to lie in some subset $\mathcal S\subset \{-1, 0, 1\}$.
	\item Theorem \ref{theorem_embeded}.1 and  Corollary \ref{corollaryPlanedAry}: labelled $\mathcal A$-ary trees, i.e. plane trees where the number of children of each vertex belongs to some fixed set $\mathcal A\subset\mathbb Z_{\geq 0}$ containing $0$ and with label increments constrained to lie in $\mathcal S\subset \{-1, 0, 1\}$. This includes in particular $d$-ary trees ($\mathcal A=\{0, d\}$) and plane trees ($\mathcal A=\mathbb Z_{\geq 0}$).
	\item Theorem \ref{theorem_Sembedded_models}.1:  canonically $\mathcal S$-embedded complete trees, with $\mathcal S\subset\{-1, 0, 1\}$. These are labelled plane trees in which every non-leaf vertex has exactly $\mathrm{Card}(\mathcal S)$ children, whose label increments, from left to right, are precisely the elements of $\mathcal S$ in increasing order.
    
	\item Theorem \ref{theorem_Sembedded_models}.2 : canonically $\mathcal S$-embedded incomplete trees, with $\mathcal S\subset \{-1, 0, 1\}$. These are labelled plane trees in which the edge increments belong to $\mathcal S$, the labels of the children of every vertex are strictly increasing from left to right, and each vertex may realize only a subset of the increments in $\mathcal S$.
	\item Theorem \ref{theorem_Sembedded_models}.3 :  canonically embedded Motzkin trees, namely labelled plane trees $(t, \ell_t)$ where each non-leaf vertex $v$ of $t$ has either one child with the same label, or exactly two children, the first one with label $\ell_t(v)-1$, and the second one with label $\ell_t(v)+1$.
\end{itemize}
More generally, the methods developed in this paper extend naturally to any class of trees that can be generated by a random tree with suitably chosen offspring and displacement distributions.

\section{Main results}

\subsection{Models of random trees.} 
\subsubsection{Well-labelled trees.} A \emph{rooted tree} $t$ is a connected, acyclic graph on a finite vertex set $V(t)$, with a distinguished vertex $\rho_t \in V(t)$ called the \emph{root} of $t$. The root induces a genealogical structure on $t$: the set of children of a vertex $v\in V(t)$ is the (possibly empty) set of vertices adjacent to $v$ that are farther from $\rho_t$ than $v$ with respect to the graph distance in $t$. We say that $t$ is a \emph{plane} tree if, in addition, for every vertex $v \in V(t)$, a linear order is specified on the set of its children. For every $v \in V(t)$, we write $c_t(v)$ for the number of children of $v$. Moreover, if $v\neq \rho_t$, we write $\pi_t(v)$ for the only vertex having $v$ as a child, that we call the parent of $v$. By convention, the edges of $t$ are oriented away from the root, that is, from parents to children (note that this convention may differ for other references). This orientation induces a partial order $\preceq$ on the vertices of $t$, where $u \preceq v$ if and only if $u$ lies on the ancestral path from the root to $v$. We write $u \prec v$ when $u \preceq v$ and $u \ne v$. In this work, we are interested in models of \emph{well-labelled} rooted plane trees, defined as follows.

\begin{definition}
Let $x\in \mathbb{Z}$. A well-labelled tree with the root labelled $x$ is a pair $(t,\ell_{t})$, where
\begin{itemize}[label=\scriptsize$\bullet$]
\item $t$ is a rooted plane tree; and
\item $\ell_{t} : V(t)\to \mathbb Z$ is a labelling function assigning an integer to each vertex of $t$, such that $\ell_{t}(\rho_{t})=x$ and for every non-root vertex $v$,
\begin{align*}
 \ell_{t}(v)-\ell_{t}\left(\pi_{t}(v)\right) \in \{-1,0,1\}.
 \end{align*}
\end{itemize}
We write $\mathbb T_x$ for the set of all well-labelled trees rooted at $x$.
\end{definition}
  If $t\in \mathbb T_x$, we will (somewhat abusively) keep the notation $t$ to talk about the associated rooted plane tree (so that the notation $V(t), \rho_t,  c_t(v), \pi_t(v)\ldots$ still make sense), and we will always write $\ell_t:V(t)\to \mathbb Z$ for the associated labelling function. Let 
  \begin{align*}
  	\Delta= \{(i, j)\in \mathbb Z^2 \mid |i-j|\leq 1\},
  \end{align*}
which is exactly the set of label pairs that can occur at the endpoints of an edge in a well-labelled tree. The main object of interest in this paper is defined as follows.
\begin{definition}
	Let $x\in \mathbb Z$ and $t\in \mathbb T_x$. 
	For every $(i, j)\in \Delta$, let $e^t_{i,j}$ denote the number of edges oriented from a vertex with label $i$ to a vertex with label $j$ in $t$, that is
	\[ e^t_{i,j} \coloneqq \mathrm{Card} \left\{
	v\in V(t)\setminus\{\rho_t\}:
	\left(\ell_t(\pi_t(v)),\ell_t(v)\right)=(i,j) \right\}.
	\]
	The family $(e^t_{i,j})_{(i,j)\in\Delta}$ is called the \emph{vertical edge profile} of $t$.
\end{definition}
 For every well-labelled tree $(t, \ell_t)$ and every $v\in V(t)$, let
\begin{align*}
	[v]_t=\Big(\ell_t(v_1)-\ell_t(v), \ell_t(v_2)-\ell_t(v), \cdots,  \ell_t(v_{c_t(v)})-\ell_t(v)\Big),
\end{align*} 
where $v_1, \cdots, v_{c_t(v)}$ are the children of $v$, listed from left to right. If $v$ is a leaf, we set $[v]_t=\varnothing$. In other words, $[v]_t$ is the vector of label increments along the edges from $v$ to its children.
\subsubsection{Random labelled trees.}
We now introduce the class of random labelled trees considered in this work. They correspond to labelled Bienaymé--Galton--Watson trees, such that the vectors of children displacement $[v]_t$ of each vertex $v$ follow fixed displacement laws that may depend on the number of children of $v$. More precisely, fix a probability distribution $\xi$ on $\mathbb Z_{\geq0}$ called the \emph{offspring distribution}, that we assume to be non-trivial (i.e. $\xi(1)<1$) and critical or subcritical, namely
\begin{align*}
	\sum_{k\geq 0} k\, \xi(k)\leq 1.
\end{align*}
Moreover, fix a family of probability distributions $\eta=(\eta^{(k)})_{k\geq 1}$ called the \emph{displacement laws}, where, for every $k\geq 1$, $\eta^{(k)}$ is a probability distribution on $\{-1, 0, 1\}^k$.
Fix $x\in \mathbb Z$. We construct a random well-labelled tree $(T, \ell_T)\in \mathbb T_x$ in three steps. 
\begin{enumerate}
\item Sample a Bienaymé--Galton--Watson tree $T$ with offspring distribution $\xi$. Since $\xi$ is critical or subcritical, the resulting tree is almost surely finite.
\item Conditionally on $T$, for every vertex $v\in V(T)$, sample a random vector $L_v\in \{-1, 0, 1\}^{c_T(v)}$ with distribution $\eta^{(c_T(v))}$, independently over the vertices.
\item The labelling function $\ell_T$ is then defined inductively by setting $\ell_T(\rho_t)=x$ and requiring that we have $[v]_T= L_v$ for every vertex $v \in V(T)$.
\end{enumerate}
The labelled tree $(T, \ell_T)$ is called a $(\xi,\eta)$-well-labelled BGW tree rooted at $x$.
We denote by $\Pi_x$ the distribution of $(T,\ell_T)$. Explicitly, for every $t\in \mathbb T_x$,
\begin{align*}
	\Pi_x(t)= \prod_{v \in V(t)}  \xi\left(c_t(v)\right) \, \eta^{(c_t(v))}\left([v]_t\right),
\end{align*}
where by convention $\eta^{(0)}(\varnothing)=1$. 
\begin{figure}[h]
		\begin{minipage}{0.32\textwidth}
		\centering
		\includegraphics[width=0.7\textwidth]{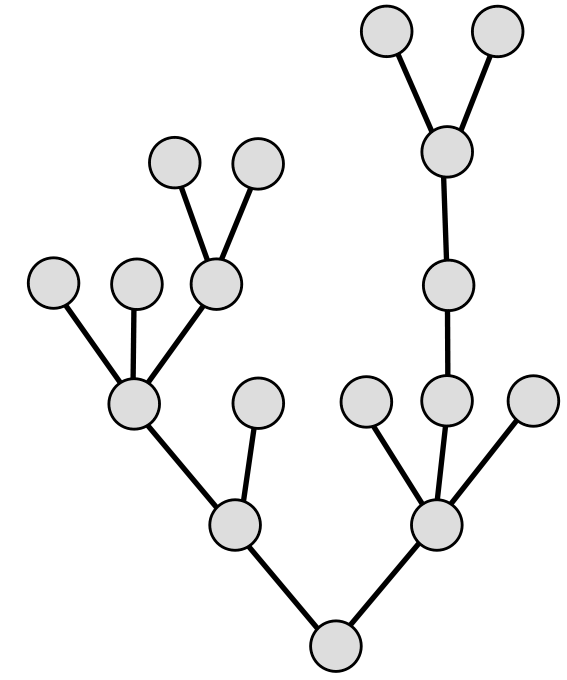}
	\end{minipage}
	\begin{minipage}{0.32\textwidth}
		\centering
		\includegraphics[width=0.7\textwidth]{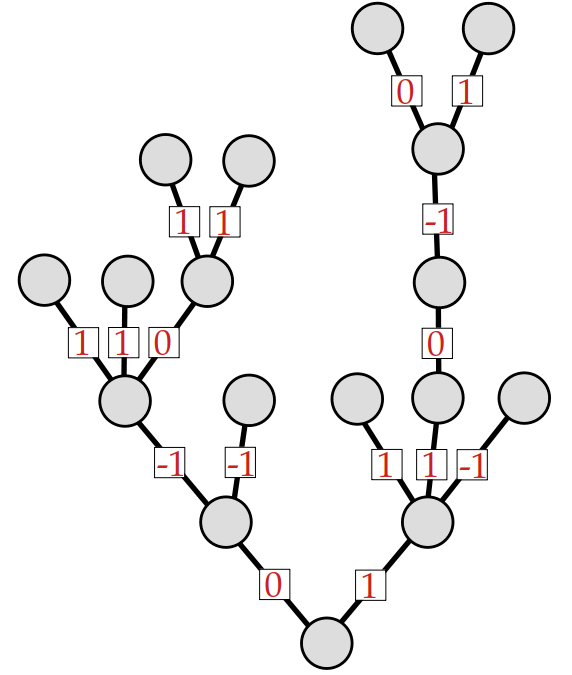}
	\end{minipage}
	\begin{minipage}{0.32\textwidth}
		\centering
		\includegraphics[width=0.7\textwidth]{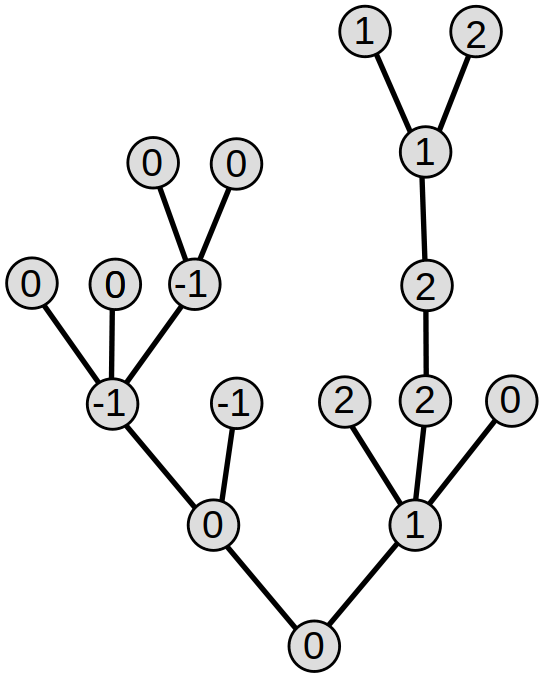}
	\end{minipage}
    \caption{Sampling of a $(\xi, \eta)$-well-labelled tree. Step $1$ (left): sample a Bienaymé--Galton--Watson tree $T$ with offspring distribution $\xi$; Step 2 (middle): independently sample the label increments along the edges from each vertex $v$ to its children according the law $\eta^{(c_T(v))}$. Step 3 (right): each vertex is labelled by the sum of the increments along the unique path from the root.}
\end{figure} 
\subsection{Markov properties for the vertical edge profile.} 
Let $T$ be a random labelled tree with law $\Pi_0$. It follows from \cite[Theorem 4.2]{MarkovPropertyMD} that
\begin{align*}
(e_{k-1, k}^T, e_{k, k-1}^T)_{k\geq 1}, \qquad \text{ and } \qquad 	(e_{-k+1, -k}^T, e_{-k, -k+1}^T)_{k\geq 1},
\end{align*}
 are time-homogeneous Markov chains which are independent conditionally on the value taken by $(e_{0, 1}^T, e_{1, 0}^T)$.
 For our purposes, it will be more convenient to introduce the (enriched) processes ${\mathcal Y}^+$ and  ${\mathcal Y}^-$ given for every $k\geq 1$ by
\begin{align}
	\label{DefX+}
	{\mathcal Y}^+_k& \coloneqq (e_{k-1, k}^T, e_{k-1, k-1}^T, e_{k, k-1}^T),\\
	\label{DefX-}	{\mathcal Y}^-_k& \coloneqq (e_{-k+1, -k}^T, e_{-k+1, -k+1}^T, e_{-k, -k+1}^T).
\end{align}
As we will see in Proposition \ref{MP}, the processes $ {\mathcal Y}^+$  and $\mathcal Y^-$ are also time-homogeneous Markov chains under $\Pi_0$, which take values in the state space:
\begin{align*}
\mathfrak S^* \coloneqq \{(p, q, r)\in \mathbb Z_{\geq 0}^3 \mid r=0 \text{ if } p=0\},
\end{align*}
since the existence of an edge from label $k$ to label $k-1$ necessarily implies the existence of at least one edge from label $k-1$ to label $k$.
\subsubsection{Kernel Formulas} \label{KernelFormSec}
 Let $ {\mathbf R}$ denote the transition kernel of ${\mathcal Y}^+$, that is, for every $(p, q, r)\in \mathfrak S^*$ such that there exists $m\geq 1$ with $\Pi_0({\mathcal Y}^+_m=(p, q, r))>0$ and every $a, b, c\geq 0$, we let
\begin{align*}
	{\mathbf R}_{p, q, r}^{a, b, c} \coloneqq \Pi_0\left(  {\mathcal Y}^+_{m+1}= (a, b, c) \mid {\mathcal Y}^+_m=(p,q, r)\right).
\end{align*}
Since the Markov chain ${\mathcal Y^+}$ is time-homogeneous, this quantity is independent of the choice of $m\geq 1$ such that $\Pi_0({\mathcal Y}^+_m=(p, q, r))>0$. Similarly, we denote by $ {\mathbf L}$ the transition kernel of ${\mathcal Y}^-$. For every $(p, q, r)\in \mathbb Z_{\geq 0}^3$ such that there exists $m\geq 1$ with $\Pi_0\left({\mathcal Y}_m^{-}=(p, q, r)\right)>0$, we let, for every $(a, b, c)\in \mathbb Z_{\geq 0}^3$
\begin{align*}
	{\mathbf L}_{p,q, r}^{a, b, c}&\coloneqq \Pi_0\left(  {\mathcal Y}_{m+1}^{-}= (a, b, c) \mid  {\mathcal Y}_m^{-}= (p, q, r)\right).
\end{align*}  Our first main result provides explicit formulas for the kernels ${\mathbf R}$ and ${\mathbf L}$, involving certain random walks related to our models of random trees. More precisely, the transition kernel can be written in terms of the following distributions associated with our models,
\begin{itemize} [label=\scriptsize$\bullet$]
	\item $\nu_+: \mathbb Z_{\geq 0}\to [0, 1]$ is the law of
	\begin{align*}
		\text{Card}\left\{v\in V(T) \, : \, \ell_T(v)=-1, \text{ and } \ell_T(u)>-1, \, \forall u\preceq v \right\},
	\end{align*}
	 the number of \emph{first hitting points of label $-1$} in $T$, when $T$ is sampled with law $\Pi_0$.
	 That is, for every $k\geq 0$, $\nu_+(k)$ is the probability that there are exactly $k$ vertices labelled $-1$ in $T$ whose ancestral line contains no other vertex labelled $-1$.
	\item $\nu_-: \mathbb Z_{\geq 0}\to [0, 1]$ is the law of
	\begin{align*}
		\text{Card}\left\{v\in V(T) \, : \, \ell_T(v)=1, \text{ and } \ell_T(u)<1, \, \forall u\preceq v \right\},
	\end{align*}
	the number of \emph{first hitting points of label $1$} in $T$.
	\item $\eta_*: (\mathbb Z_{\geq 0})^3 \to [0, 1]$ is the distribution of the numbers of children of the root labelled $-1, 0$ and $1$. That is, for every $(a, b, c)\in (\mathbb Z_{\geq 0})^3$ the value of $\eta_*(a, b, c)$ is the probability under $\Pi_0$ that exactly $a$ children of the root have label $-1$, $b$ children of the root have label $0$ and $c$ children of the root have label $1$.
\end{itemize}
The formula we derive make use of convolutions of these laws.  For every $n, m\geq 0$, we let
\begin{align}
	f^{\nu_+}_n(m)= \mathbb P\left(\sum_{i=1}^n N_i^+= m\right) \qquad  \text{ and } \qquad f^{\nu_-}_n(m)= \mathbb P\left(\sum_{i=1}^n N_i^-= m\right),
\end{align} where $(N_i^+)_{i\geq 1}$ (resp. $(N_i^-)_{i\geq 1}$) is a sequence of i.i.d. random variables with law $\nu_+$ (resp. $\nu_-$). Similarly, for every $n, x, y , z \in \mathbb Z_{\geq 0}$, let
\begin{align}
	f^{\eta_*}_n(x, y , z)=	\mathbb P\left(\sum_{i=1}^n (x_{i}^{(-1)},x_i^{(0)}, x_i^{(1)})=(x, y, z)\right),
\end{align}
with $(x^{(-1)}_i, x^{(0)}_i, x_{i}^{(1)})_{i\geq 1}$ a sequence of i.i.d. random variables with law $\eta_*$. The kernels admit the following explicit expressions.
	\begin{theorem}[Kernel formulas]\label{FormKernel} \, \medskip 
		\begin{enumerate}
			\item For every $(p,q, r)\in \mathbb Z^3_{\geq 0}$ such that there exists $m\geq 1$ with $\Pi_0\left( {\mathcal Y}_m^{+}=(p, q, r)\right)>0$, we have $f_p^{\nu_+}(r)>0$. Moreover, for every $(a, b, c)\in \mathbb Z^3_{\geq 0}$,
			\begin{align}
				\label{FormKernelEq}  {\mathbf R}^{a, b, c}_{p,q, r}= \frac{p}{p+b+c}f^{\eta_*}_{p+b+c}(r, b, a)\frac{f^{\nu_+}_{a}(c)}{f^{\nu_+}_{p}(r)},
			\end{align}
			where, by convention, ${\mathbf R}_{0, q, 0}^{a, b, c}=\1_{(a, b, c)=(0, 0, 0)}$, for every $q\geq 0$.
			\item For every $(p,q, r)\in \mathbb Z_{\geq 0}^3$ such that there exists $m\geq 1$ with $\Pi_0\left( {\mathcal Y}_m^{-}=(p, q, r)\right)>0$, we have $f_p^{\nu_-}(r)>0$. Moreover, for every $(a, b, c)\in \mathbb Z^3_{\geq 0}$,
			\begin{align}
				\label{FormKernelEq2}
				{\mathbf L}^{a,b,c}_{p,q,r}= \frac{p}{p+b+c}f^{\eta_*}_{p+b+c}(a, b, r)\frac{f^{\nu_-}_{a}(c)}{f^{\nu_-}_{p}(r)},
			\end{align}
			where, by convention, ${\mathbf L}_{0, q, 0}^{a, b, c}=\1_{(a, b, c)=(0, 0, 0)}$ for every $q\geq 0$.
		\end{enumerate} 
\end{theorem}

\begin{remark} 	In many classical models, the factor $f^{\eta_*}_{p+b+c}(r, b, a)$ is completely explicit and the kernel formula simplifies nicely. For instance:
	\begin{enumerate}
		\item The formula \eqref{FormKernelEq} generalizes the kernel formula obtained in \cite[Theorem 5.1]{MarkovPropertyMD} for the canonically embedded incomplete binary tree model, which corresponds to the offspring distribution $\xi=\frac{1}{4}\delta_{0}+\frac{1}{2}\delta_1+\frac{1}{4}\delta_2$ and the displacement laws $\eta^{(1)}= \frac{1}{2}\delta_{-1}+\frac{1}{2}\delta_1$ and $\eta^{(2)}= \delta_{(-1, 1)}$. In this setting, note that we must have $b=q=0$ for the transition kernel \eqref{FormKernelEq} to make sense and be non-zero (since necessarily $e^T_{k, k}=0$ for every $k\in \mathbb Z$). Moreover, a straightforward computation yields $f^{\eta_*}_{n}(i, 0, j)= 4^{-n}{n\choose i}{n\choose j}$. Substituting this into Equation \eqref{FormKernelEq}, we get
		\begin{align*}
			{\mathbf R}^{a, 0, c}_{p,0, r}=  {\mathbf L}^{a,0,c}_{p,0,r}= \frac{p4^{-p-c}}{p+c}{p+c\choose r}{p+c \choose a}\frac{f^{\nu_+}_{a}(c)}{f^{\nu_+}_{p}(r)}.
		\end{align*}
		\item When the displacements are uniform in $\{-1,1\}$ (that is $\eta^{(d)}=\mathcal U(\{-1, 1\})^{\otimes d})$, then again we must have $b=q=0$ for the kernel to make sense and be non-zero. We easily check that $f^{\eta_*}_{n}(i, 0, j)= f^\xi_n(i+j)2^{-i-j}{i+j\choose i}$ and as a direct consequence of Equation \eqref{FormKernelEq} we have
		\begin{align*}
			{\mathbf R}^{a, 0, c}_{p,0, r}= {\mathbf L}^{a,0,c}_{p,0,r}=\frac{p2^{-r-a}}{p+c}f^{\xi}_{p+c}(r+a){r+a\choose a}\frac{f^{\nu_+}_{a}(c)}{f^{\nu_+}_{p}(r)},
		\end{align*}
		(an analogous formula holds when the displacements are uniform in $\{-1, 0, 1\}$). 
	\end{enumerate} 	
\end{remark}
\subsection{Initial laws and a path probability formula} We also derive explicit formulas for the laws of $(e_{1, 0}^T,e_{0, 1}^T)$ and $(\mathcal Y^+_1, \mathcal Y^-_1)$.

\begin{theorem} 
	\label{th_profat0} For every $p, r\in \mathbb Z_{\geq 0}$, we have: 
	\begin{align} 
		\label{LawiniSemi}
		\Pi_0\left(e_{0,1}^T= p , e_{1,0}^T= r\right)= \frac{f^{\nu_-}_{r+1}(p)f^{\nu_+}_{p}(r)}{r+1}.
	\end{align}
	Moreover, for every $p, q, r, a,b, c\in \mathbb Z_{\geq 0}$, we have:
	\begin{align}
		\label{LawinitialCondi}
		\Pi_0\left({\mathcal Y}_1^+=(p, q, r), {\mathcal Y}_1^{-}= (a, b, c)\right)= \frac{f^{\nu_-}_{a}(c)f^{\nu_+}_{p}(r)}{c+b+r+1}f^{\eta_*}_{c+b+r+1}(a,b,p)\mathds 1_{b=q}.
	\end{align}
\end{theorem}

\begin{remark}
	We can adapt the proof of Theorem \ref{th_profat0} for models of forests. More precisely, for every $k\geq 1$, let $F_k$ be a forest of $k$ independent labelled trees with law $\Pi_0$ and write $\Pi_0^{\otimes k}$ for its law. For every $(i, j)\in \Delta$, we write $e^{F_k}_{i, j}$ for the total number of edges labelled $i\rightarrow j$ in $F_k$. The proof of Theorem \ref{th_profat0} extends directly to forests, yielding for every $p, r\in \mathbb Z_{\geq 0}$:
	\begin{align*}
		\Pi_0^{\otimes k}\left(e_{0,1}^{F_k}= p , \ e_{1,0}^{F_k}= r \right)= \frac{k\, f^{\nu_-}_{r+k}(p)\, f^{ \nu_+}_{p}(r)}{r+k}.
	\end{align*}
\end{remark}
Using the formulas provided by Theorem \ref{FormKernel} and Theorem \ref{th_profat0}, we can derive a general product-form expression for the probability that a tree has a given vertical edge profile. Remarkably, the resulting expression no longer involves factors using $f^{\nu_-}$ and $f^{\nu_+}$.
\begin{proposition}[Path probabilities]
	\label{Path_probabilities}
	Let $(x_{i, j})_{(i, j)\in \Delta}$ be a collection of non-negative integers with finitely many non-zero elements.  For every $k\in \mathbb Z$, write:
	\[n_k= \1_{k= 0}+ x_{k-1, k}+x_{k, k}+ x_{k+1, k},\]
	and assume that $n_i>0$ for every $(i, j)\in \Delta$ such that $x_{i, j}>0$. Let
	\begin{align*}
		M^*= \max\{k\geq 0 \mid n_k>0\}, \quad \text{ and } \quad M_*=\min\{k\leq 0 \mid n_k>0\}.
	\end{align*}
	and assume that $n_k>0$ for every $M_*\leq k\leq M^*$.  We have
	\begin{align}
		\label{FormulPaths}
		\Pi_0\left(e_{i, j}^T=x_{i, j}, \ \forall (i, j)\in \Delta\right)= \prod_{k=M_*}^{M^*} \left[\frac{1}{n_k}f^{\eta_*}_{n_k}(x_{k, k-1}, x_{k, k}, x_{k,k+1})\right] \prod_{k=1}^{M^*}x_{k-1, k} \prod_{k=M_*}^{-1} x_{k+1,k}, \qquad
	\end{align}
	where the last two products are understood to equal $1$ when the index set is empty.
\end{proposition}

\begin{remark}
	The numbers $n_k$ in the statement of Proposition \ref{Path_probabilities} simply correspond  to the number of vertices labelled $k$ in a tree with vertical edge profile $(x_{i, j})_{(i, j)\in \Delta}$.
\end{remark}

\subsection{Probabilistic counting of trees}
\label{ProbaCounting}
The path probability formula of Proposition \ref{Path_probabilities} can be used as an enumeration tool. The idea is that for suitable choices of offspring and displacement laws, the associated random tree conditioned on its size is uniformly distributed over interesting combinatorial classes of labelled trees. In that case, the probability to have a given vertical edge profile translates directly into formulas for the number of trees with a prescribed vertical edge profile in the associated combinatorial class. We illustrate this principle on several families of embedded trees.

\subsubsection{$\mathcal S$-embedded Cayley and $\mathcal A$-ary trees}
\label{ProbaCountingSembedded}
Our first results concern Cayley and $\mathcal A$-ary trees, defined as follows.
\begin{itemize}[label=\scriptsize$\bullet$]
	\item A rooted Cayley tree (or unordered tree) $t$ of size $n$ is a rooted tree whose vertex set is $[n]=\{1, \cdots, n\}$,
	\item Fix $\mathcal A\subset \mathbb Z_{\geq0}$ a set containing $0$. An $\mathcal A$-ary tree of size $n$ is a plane tree $t$ of size $n$ whose vertices all have out-degree in $\mathcal A$, that is:
	\begin{align*}
		\forall v\in V(t), \qquad c_t(v)\in  \mathcal A.
	\end{align*}
	When $\mathcal A=\{0, d\}$ for some $d\geq 1$, we simply say that $t$ is a $d$-ary tree. Moreover a $\mathbb Z_{\geq0}$-ary tree is simply a plane tree with no additional constraint.
\end{itemize}
For every integer $n, k\geq 0$, write $S_{\mathcal A}(n, k)$ for the number of ways to write $k$ as an (ordered) sum of $n$ elements of $\mathcal A$. Equivalently:
\begin{align*}
	S_{\mathcal A}(n, k)= [z^k]\left(\sum_{a\in \mathcal A}z^a\right)^n.
\end{align*}
In particular for every $d\geq 1$ and every $n, k\geq 0$:
\begin{align}
	\label{ExplicitSAdary} S_{\{0, d\}}(n , k)= \left\{ \begin{array}{l} {n\choose k/d} \quad \text{ if } d|k, \ n\geq k/d, \\
		0 \qquad \ \ \text{ otherwise}, \end{array} \right.
\end{align} and for every $n, k\geq 0$
\begin{align}
	\label{ExplicitSplane}
	S_{\mathbb Z_{\geq 0}}(n , k)={n+k-1  \choose n-1}.
\end{align}
\begin{definition}
	Fix a non-empty subset $\mathcal S\subset\{-1, 0, 1\}$. An $\mathcal S$-embedded rooted Cayley tree (resp. plane tree, resp. $d$-ary, resp. $\mathcal A$-ary tree) is a pair $(t, \ell_t)$ where:
	\begin{itemize}[label=\scriptsize$\bullet$]
		\item $t$ is a rooted Cayley tree (resp. a plane tree,  resp. $d$-ary, resp. an $\mathcal A$-ary tree),
		\item $\ell_t: V(t)\to \mathbb Z$ is a labelling such that $\ell_t(\rho_t)=0$ and for every $v\in V(t)\setminus\{\rho_t\}$:
		\begin{align*}
			\ell_t(v)-\ell_t(\pi_t(v))\in \mathcal S.
		\end{align*}
	\end{itemize}
\end{definition}
\noindent Fix $\mathcal S\subset\{-1, 0, 1\}$ and a class of trees $\mathcal C\in \{\text{Cayley}, \mathcal A\text{-ary},\text{plane}, d\text{-ary}\}$. Let $\underline x= (x_{i, j})_{(i, j)\in \Delta_{\mathcal S}}$ be a collection of non-negative integers with finitely many non-zero elements, indexed by:
\begin{align*}
	\Delta_{\mathcal S}=\{(i, j)\in \mathbb Z \mid j-i\in \mathcal S\}.
\end{align*}
An $\mathcal S$-embedded $\mathcal C$-tree $(t, \ell_t)$ is said to have vertical edge profile $\underline x$ if for every $(i, j)\in \Delta_{\mathcal S}$, the number of edges in $t$ oriented from a vertex labelled $i$ to a vertex labelled $j$ is exactly $x_{i, j}$ (where edges are again oriented away from the root). We write $\mathcal T_{\mathcal C}^{\mathcal S}(\underline x)$ for the set of all $\mathcal S$-embedded $\mathcal C$-trees with vertical edge profile $\underline x$.
The following theorem gives enumeration formulas for its cardinality.
\begin{theorem}
	\label{theorem_embeded}
	Fix $\mathcal S\subset \{-1, 0, 1\}$ and a collection $\underline x= (x_{i, j})_{(i, j)\in \Delta_{\mathcal S}}$ of non-negative integers with finitely many non-zero elements. Set $x_{i, j}=0$ whenever $j-i\notin \mathcal S$.  For every $k\in \mathbb Z$, let 
	\begin{align*}
		n_k= \1_{k= 0}+ \sum_{s\in \mathcal S} x_{k-s, k},  \qquad \text{and} \qquad  \kappa_k= \sum_{s\in \mathcal S} x_{k, k+s},
	\end{align*}
	and write $M^*= \max\{k\geq 0 \mid n_k>0\}$ and $M_*=\min\{k\leq 0 \mid n_k>0\}$. Assume that $n_k>0$ for every $M_*\leq k\leq M^*$ and that $n_i>0$ for every $(i, j)\in \Delta$ such that $x_{i, j}>0$
	\begin{enumerate}
		\item Let $\mathcal A\subset \mathbb Z_{\geq 0}$ with $0\in \mathcal A$ and $\mathcal A\neq \{0\}$. We have:
		\begin{align}
			\label{EnumAary}
			\mathrm{Card}\bigg(\mathcal T^{\mathcal S}_{\mathcal A\text{-ary}}\left(\underline x\right)\bigg)=  \frac{\displaystyle \prod_{1\leq k\leq M^*}x_{k-1,k} \prod_{M_*\leq k\leq -1}x_{k+1,k}}{\displaystyle \prod_{(i, j)\in \Delta_{\mathcal S}} x_{i, j}!} \prod_{M_*\leq k\leq M^*}\frac{\kappa_k!}{n_k}S_{\mathcal A}(n_k, \kappa_k),
		\end{align}
		\item Let $N=\sum_{i, j}x_{i, j}$ be the total number of edges. Then
		\begin{align}
			\label{EnumCayley}
			\mathrm{Card}\bigg(\mathcal T^{\mathcal S}_{\text{Cayley}}\left(\underline x\right)\bigg)=  \frac{\displaystyle \prod_{1\leq k\leq M^*}x_{k-1,k} \prod_{M_*\leq k\leq -1}x_{k+1,k}}{\displaystyle \prod_{(i, j)\in \Delta_{\mathcal S}} x_{i, j}!}\left(N+1\right)!& \prod_{M_*\leq k\leq M^*} n_k^{\kappa_k-1}.
		\end{align}
	\end{enumerate}
\end{theorem}
\begin{remark}
	Formula \eqref{EnumCayley} was already proven in \cite[Theorem $5$]{BousquetChapuy}, via bijective methods, and reproved in \cite{Bernardi} through non bijective but combinatorial methods.
\end{remark}
The quantities $n_k$ and $\kappa_k$ admit a simple combinatorial interpretation. If a plane tree $t$ has vertical edge profile $\underline x$, then $n_k$ is the number of vertices labelled $k$ and $\kappa_k$ is the number of non-root vertices with parent labelled $k$.
\begin{figure}[h!]
	\centering\includegraphics[width=0.6\textwidth]{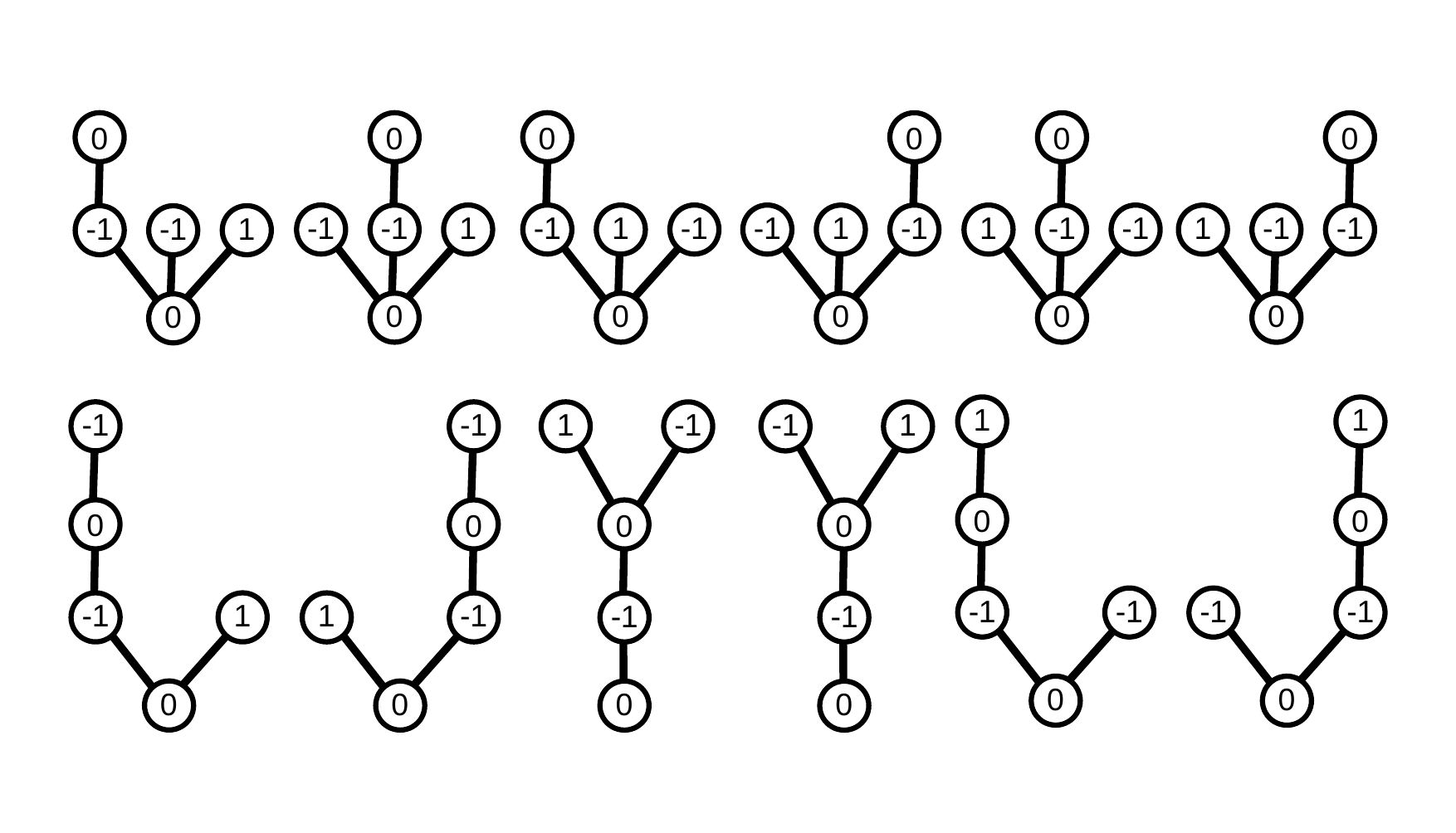}
	\caption{The twelve $\{-1, 1\}$-embedded plane trees with vertical edge profile given by $e_{0, 1}^t=1$, $e_{0, -1}^t=2$, $e_{-1, 0}^t=1$ and $e_{i, j}^t=0$ for every other $i, j$.}
\end{figure}
\begin{corollary}
	\label{corollaryPlanedAry}
	In the notation of Theorem \ref{theorem_embeded}.
	\begin{enumerate}
		\item The number of $\mathcal S$-embedded plane trees with vertical edge profile $\underline x=(x_{i, j})_{(i,j)\in \Delta_{\mathcal S}}$ is:
		\begin{align}
			\label{EnumPlane}
			\mathrm{Card}\left(\mathcal T^{\mathcal S}_{\text{plane}}\left(\underline x\right)\right)=  \frac{\displaystyle \prod_{1\leq k\leq M^*}x_{k-1,k} \prod_{M_*\leq k\leq -1}x_{k+1,k}}{\displaystyle \prod_{(i, j)\in \Delta_{\mathcal S}} x_{i, j}!} \prod_{M_*\leq k\leq M^*} \frac{(\kappa_k+ n_k-1)!}{n_k!}.
		\end{align}
		\item 	The number of $\mathcal S$-embedded $d$-ary trees with vertical edge profile $\underline x$ is:
		\begin{align}
			\label{Enumdary}
			\mathrm{Card}\left(\mathcal T^{\mathcal S}_{d\text{-ary}}\left(\underline x\right)\right)=  \frac{\displaystyle \prod_{1\leq k\leq M^*}x_{k-1,k} \prod_{M_*\leq k\leq -1}x_{k+1,k}}{\displaystyle \prod_{(i, j)\in \Delta_{\mathcal S}} x_{i, j}!} \prod_{M_*\leq k\leq M^*}\frac{\kappa_k!}{n_k}{n_k \choose \kappa_k/d},
		\end{align}
		if $d|\kappa_k$ for every $k\in \mathbb Z$ and $\mathcal T^{\mathcal S}_{d\text{-ary}}\left(\underline x\right)$ is empty otherwise.
	\end{enumerate}
\end{corollary}

\begin{remark}
	These results have a few interesting specializations:
	\begin{enumerate} 	
		\item If $\mathcal S=\{0\}$ then the vertical edge profile of a $\{0\}$-embedded tree can only be of the form $x_{0, 0}=m$ for some $m\geq 0$ and $x_{i, j}=0$ for any other $(i, j)\in \Delta$. The formulas then simply recover the number of trees of size $m$ in the class considered. That is \eqref{EnumPlane} simplifies to $\frac{1}{m+1}{2m \choose m}$, \eqref{Enumdary} simplifies to $\frac{1}{m+1}{m+1 \choose m/d}$ and \eqref{EnumCayley} simplifies to $(m+1)^{m}$.
		\item If $\mathcal S=\{1\}$ then the vertical profile coincides with the horizontal profile of the tree. That is, if $\underline x$  is a vertical edge profile of an $\{1\}$-embedded tree, then $x_{i, j}= 0$ whenever $i<0$ or $j\neq i+1$, and $x_{i, i+1}$ counts exactly the number of vertices at distance $i+1$ from the root. The formulas then give the number of trees in the class considered with prescribed number of vertices at each height (often called the horizontal profile of the tree).
	\end{enumerate}
\end{remark}

\begin{remark}
	We mention that these results can be extended to the enumeration of forests of $f$ trees with prescribed vertical edge profile: just replace $n_0$ by 
	\begin{align*}
		n_0'= f+ \sum_{s\in \mathcal S}x_{-s, 0},
	\end{align*} in the enumeration formulas.
\end{remark}

\subsubsection{Canonically embedded trees}	
\label{CanonicallyEmbeddedSec}
\begin{figure}[h]
	\begin{minipage}{0.32\textwidth}
		\centering	\includegraphics[width=0.6\textwidth]{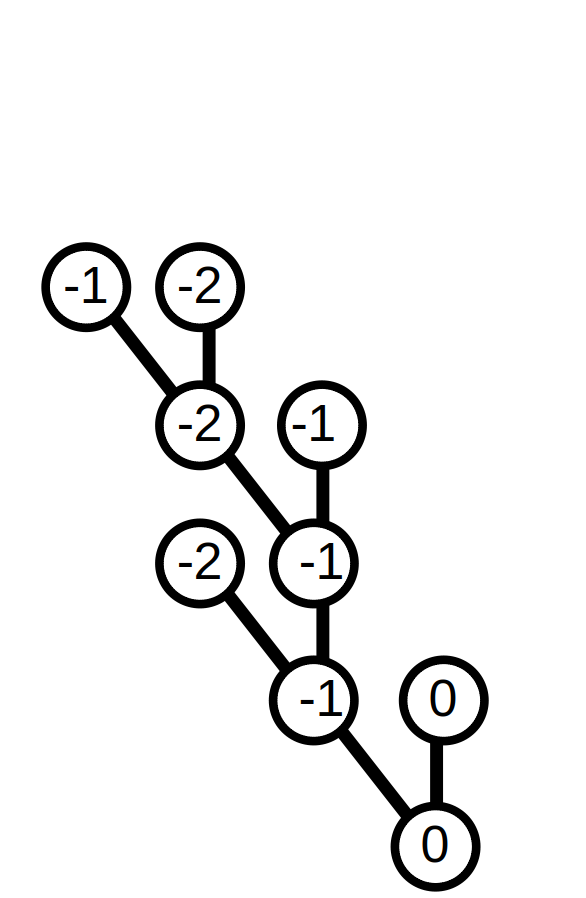}
	\end{minipage}\vrule
	\begin{minipage}{0.32\textwidth}
		\centering\includegraphics[width=0.85\textwidth]{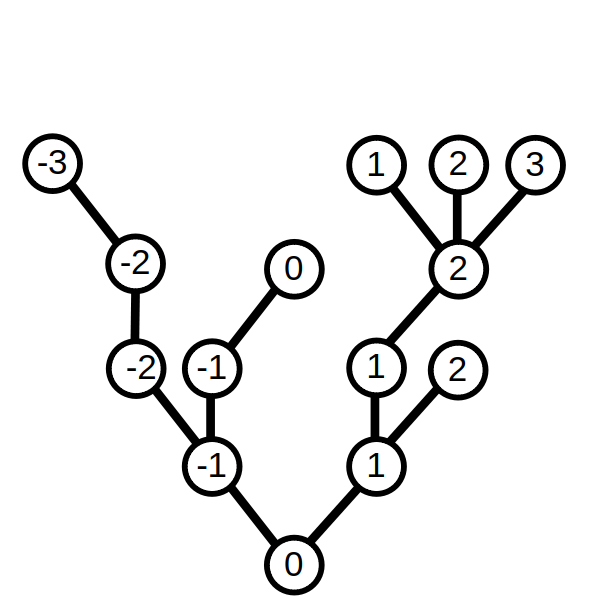}
	\end{minipage}\vrule
	\begin{minipage}{0.32\textwidth}
		\centering	\includegraphics[width=0.6\textwidth]{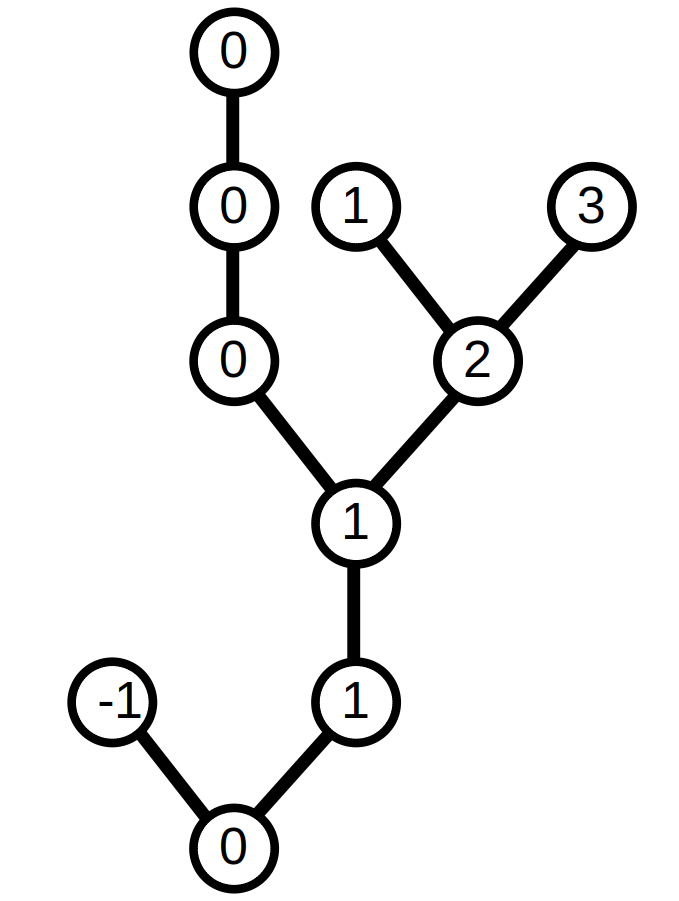}
	\end{minipage}
	\caption{A canonically $\{-1, 0\}$-embedded complete tree (left), a canonically $\{-1, 0, 1\}$-embedded incomplete tree (middle) and a canonically embedded Motzkin tree (right).}
	\label{FigExamEmbb}
\end{figure}
We now describe several models of labelled trees where the labelling is determined by the structure of the tree.
Let $\mathcal S\subset\{-1, 0, 1\}$ with cardinal $|\mathcal S|\geq 2$ and fix $\underline x= (x_{i, j})_{(i, j)\in \Delta_{\mathcal S}}$ a collection of non-negative integers with finitely many nonzero terms.
\paragraph{Canonically $\mathcal S$-embedded complete trees} A canonically $\mathcal S$-embedded complete tree is a $|\mathcal S|$-ary tree $(t, \ell_t)$ such that for every non-leaf vertex $v\in V(t)$, we have
\begin{align*}
	[v]_t=\mathcal S,
\end{align*} where we recall that $[v]_t$ is the vector of label increments from $v$ to its children, in left-to-right order. For instance, if $\mathcal S= \{-1, 1\}$ then the tree is a complete binary tree and the first child (resp. the second child) of each vertex that is not a leaf must have label $1$ less (resp. $1$ more) than its parent.  We write $\mathcal T^{\mathcal S}_{\text{comp.}}(\underline x)$ for the set of all canonically $\mathcal S$-embedded complete trees with vertical edge profile $\underline x$. 

\paragraph{Canonically $\mathcal S$-embedded incomplete trees}
An incomplete canonically $\mathcal S$-embedded tree is an $\mathcal S$-embedded plane tree $(t, \ell_t)$ such that, for every vertex $v\in V(t)$, the labels of its children are strictly increasing in left-to-right order. We write $\mathcal  T^{\mathcal S}_{\text{incomp.}}(\underline x)$ for the set of all incomplete canonically $\mathcal S$-embedded trees with vertical edge profile $\underline x$.

\paragraph{Canonically embedded Motzkin trees} A canonically embedded Motzkin tree is a labelled plane tree $(t, \ell_t)$ where each vertex has at most $2$ children. Moreover for every vertex $v\in V(t)$, if $v$ has exactly one child $c$ then $\ell_t(c)= \ell_t(v)$ and if $v$ has two children $(c_1, c_2)$ (in the order prescribed by the plane structure) then $\ell_t(c_1)=\ell_t(v)-1$ and $\ell_t(c_2)= \ell_t(v)+1$.  If $\underline x= (x_{i, j})_{(i, j)\in \Delta}$ is a collection of non-negative integers, we write $\mathcal T_{\text{Motzkin}}(\underline x)$ for the set of all canonically embedded Motzkin trees with vertical edge profile $\underline x$.

\begin{theorem}
	\label{theorem_Sembedded_models}
	Fix a subset $\mathcal S\subset \{-1, 0, 1\}$ with $|\mathcal S|\geq 2$ and a collection $\underline x= (x_{i, j})_{(i, j)\in \Delta_{\mathcal S}}$ of non-negative integers with finitely many non-zero elements. Set $x_{i, j}=0$ whenever $(i, j)\notin \Delta_{\mathcal S}$. For every $k\in \mathbb Z$, let 
	\begin{align*}
		n_k= \1_{k= 0}+ \sum_{s\in \mathcal S} x_{k-s, k} \qquad \text{ and } \qquad  \kappa_k= \sum_{s\in \mathcal S} x_{k, k+s}.
	\end{align*} Let $M^*= \max\{k\geq 0 \mid n_k>0\}$, and $M_*=\min\{k\leq 0 \mid n_k>0\}$.  Assume that $n_k>0$ for every $M_*\leq k\leq M^*$ and that $n_i>0$ for every $(i, j)\in \Delta$ with $x_{i, j}>0$.
	\begin{enumerate}
		\item If for every $k\in \mathbb Z$, $x_{k, k+s}$ does not depend on $s\in \mathcal S$ and equals some $y_k\in \mathbb Z_{\geq0}$, then
		\begin{align}
			\label{EnumComp}
			\mathrm{Card}\left(\mathcal T^{\mathcal S}_{\text{comp.}}(\underline x)\right)= 
			\prod_{k=1}^{M^*}x_{k-1,k}\prod_{k=M_*}^{-1}x_{k+1,k}\prod_{k=M_*}^{M^*}\frac{1}{n_k}{n_k\choose y_k},
		\end{align}
		otherwise $\mathcal T^{\mathcal S}_{\text{comp.}}(\underline x)$ is empty.
		\item We have:
		\begin{align}
			\label{EnumIncomp}
			\mathrm{Card}\left(\mathcal T^{\mathcal S}_{\text{incomp.}}(\underline x)\right)= 
			\prod_{k=1}^{M^*}x_{k-1,k}\prod_{k=M_*}^{-1}x_{k+1,k}\prod_{k=M_*}^{M^*}\frac{1}{n_k}\prod_{(i, j)\in \Delta_{\mathcal S}}{n_i\choose x_{i, j}}.
		\end{align}
		\item  If $\mathcal S=\{-1, 0, 1\}$ and if $x_{k, k-1}= x_{k, k+1}$ and $x_{k, k}+x_{k, k+1}\leq n_k$ for every $k\in \mathbb Z$, we have:
		\begin{align}
			\label{EnumMotzkin}
			\mathrm{Card}\left(\mathcal T_{\text{Motzkin}}(\underline x)\right)= 
			\prod_{k=1}^{M^*}x_{k-1,k}\prod_{k=M_*}^{-1}x_{k+1,k}\prod_{k=M_*}^{M^*}\frac{(n_k-1)!}{x_{k, k+1}!x_{k, k}!(n_k-x_{k, k}-x_{k, k-1})!}. \, \, \, \, \, \, \, \,\, \, \, \,
		\end{align}
		Otherwise $\mathcal T_{\text{Motzkin}}(\underline x)$ is empty.
	\end{enumerate}
\end{theorem}

\begin{figure}[h]
	\centering	\includegraphics[width=0.6\textwidth]{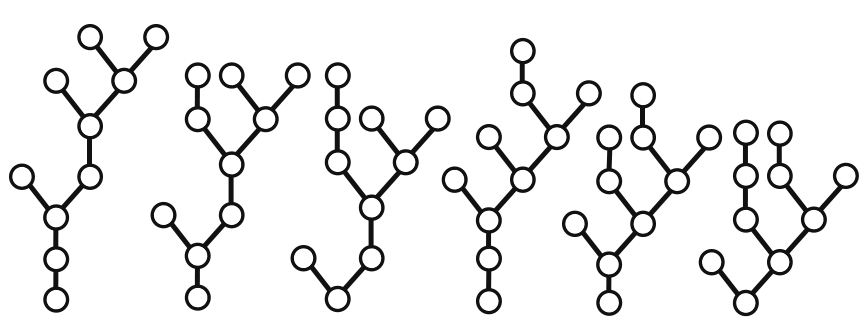}
	\caption{The six Motzkin trees with the same vertical edge profile as the example on the right-hand side of Figure \ref{FigExamEmbb}.}
\end{figure}
\begin{remark}
	Formula \eqref{EnumIncomp} was already proven in \cite[Theorem $6$]{BousquetChapuy} through bijective methods. We also mention that formula \eqref{EnumComp} can be derived directly from \cite[Theorem $4$]{BousquetChapuy} which gives a product formula for the number of canonically $\mathcal S$-embedded incomplete trees with a given vertical \emph{vertex} profile (that is, such that the number of vertices labelled $k$ is prescribed for every $k\in \mathbb Z$). 
	The idea goes as follows. If $t$ is a canonically $\mathcal S$-embedded complete tree with at least one edge, we obtain a canonically $\mathcal S$-embedded \emph{incomplete} tree $\Phi(t)$ by removing all the leaves of $t$. For every $N\geq 1$ this "pruning" operation yields a bijection between the set of all $\mathcal S$-embedded complete trees with $|\mathcal S|\times N$ edges and the set of all $\mathcal S$-embedded incomplete trees with $N$ vertices. Moreover, one readily checks that the common value of $e^{t}_{k, k+s}$ for $s\in \mathcal S$ is exactly the number of vertices labelled $k$ in $\Phi(t)$. Finally, through the bijection $\Phi$, counting $\mathcal S$-embedded complete trees with a prescribed vertical edge profile is equivalent to counting $\mathcal S$-embedded incomplete trees with a prescribed vertical vertex profile. Therefore \cite[Theorem $4$]{BousquetChapuy} immediately yields \eqref{EnumComp}. Details are omitted.
\end{remark}
\begin{remark}
These formulas also extend to the enumeration of forests of $f$ trees with prescribed vertical edge profile, just replace $n_0$ by 
	\begin{align*}
		n_0'= f+ \sum_{s\in \mathcal S}x_{-s, 0}.
	\end{align*}
\end{remark}
\section{Preliminaries}

\subsection{The excursion forest decomposition}

We first present a tree transformation that will play a central role throughout the paper.This decomposition was introduced in \cite{MarkovPropertyMD} and will again play a central role here. We follow the notation and terminology of that paper.
\begin{definition}
	Let $x\in \mathbb Z$ and $t\in \mathbb T_x$, and fix $l\in \mathbb Z$. The 
	\emph{truncation} of $t$ at label $l$, denoted by $t^{[l]}$, is obtained from $t$ by deleting all vertices of $t$ with a strict ancestor of label $l$. Equivalently, $t^{[l]}$ is the subtree of $t$ whose vertex set is:
	\begin{align*}
		\Big\{u\in V(t) \;\mid\; \forall v\prec u,\ \ell_t(v)\neq l \Big\}.
	\end{align*}
	We view $t^{[l]}$ as a labelled tree (where vertices obviously keep the same labels as in $t$).    
\end{definition}
	We next define \emph{tree excursions}, which correspond to truncations at label $0$ of trees rooted at $1$ or $-1$, and which play an important role in this paper.
	\begin{definition}
		A \emph{positive tree excursion} (resp. \emph{negative tree excursion}) is a labelled plane tree whose root has label $1$ (resp. label $-1$), in which all labels are nonnegative (resp. nonpositive), and such that each vertex with label $0$ is a leaf. We denote by $\mathcal E_+$ (resp. $\mathcal E_-$) the set of all positive (resp. negative) tree excursions.
		If $\tau$ is a tree excursion, we define 
		\begin{align*}
			n_\tau = \mathrm{Card}\{v \in V(\tau) \mid \ell_\tau(v) = 0\},
		\end{align*}
		the number of vertices with label $0$ in $\tau$.
	\end{definition}
		
	\noindent For every $x\in \mathbb Z$ and $l\neq x$, we denote by $\Pi_x^{[l]}$ the law of $T^{[l]}$, the truncation of $T$ at label $l$, when $T$ is sampled according to $\Pi_x$. 
	To simplify notation, we will write $\Pi^+ \coloneqq\Pi^{[0]}_1$ and $\Pi^- \coloneqq \Pi^{[0]}_{-1}$. Note that the support of the measure $\Pi^+$ (resp. $\Pi^-$) is contained in $\mathcal E_+$ (resp. $\mathcal E_-$). Moreover, a straightforward computation shows that if $\tau\in \mathcal E_+$ is a positive tree excursion, then
	\begin{align}
		\Pi^+(\tau) = \prod_{\substack{v \in V(\tau) \\ \ell_\tau(v) \neq 0}} \xi\left(c_{\tau}(v)\right) \, \eta^{(c_\tau(v))}\left([v]_\tau\right) ,
	\end{align}
	and the analogous expression holds for $\Pi^-$ when $\tau\in \mathcal E^-$.
	
	\medskip
	
		Let $T$ be sampled with law $\Pi_1$. Observe that $n_{T^{[0]}}$ is exactly the number of first hitting points of label $0$ in $T$. In particular, it also corresponds to the number of first hitting points of label $-1$ in the labelled tree we get by translating every label of $T$ by $-1$, which has law $\Pi_0$. Recalling the definition of $\nu_+$ from Section \ref{KernelFormSec}), it follows that for every $k\geq 0$,
	\begin{align}
		\nu_+(k)= \Pi^+(n_\tau= k),
	\end{align}
	and similarly $\nu_-(k)= \Pi^-(n_\tau= k)$.
	
	\medskip
	
\begin{figure}[h] 
		\centering
		\includegraphics[scale = 0.9]{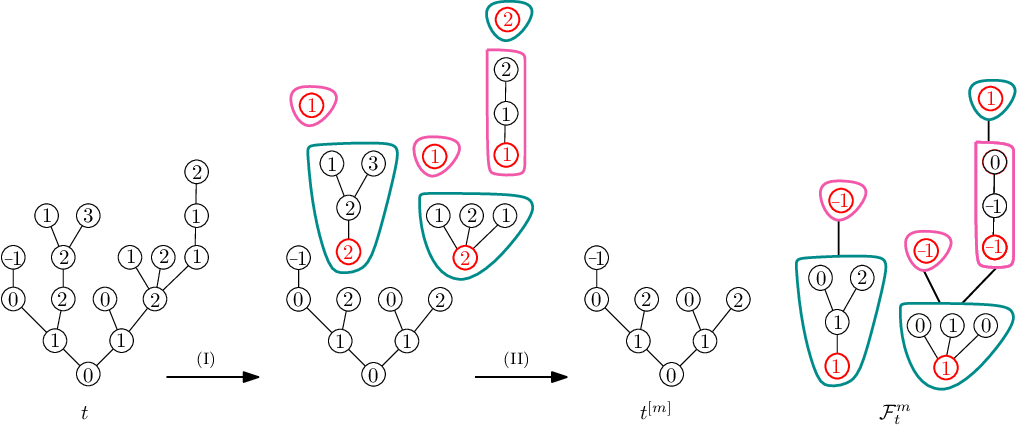}
		\caption{An example of the excursion decomposition of a tree at level $m=2$.}
		\label{Fig1}
	\end{figure}

	Let $t\in \mathbb T_0$ be a labelled plane tree rooted at $0$, and fix an integer $m \geq 1$. We next describe a decomposition of the tree $t$ into tree excursions above and below label $m$. This decomposition is illustrated in Figure \ref{Fig1}. 
	
	The first step consists in disconnecting $t$ at the endpoint of each edge that connects a vertex labelled $m-1$ to one labelled $m$, or vice versa (step (I) in Figure \ref{Fig1}). More precisely, consider an edge $e$ and denote its endpoints $u \prec v$, where $u$ is the parent of $v$. If $(\ell_t(u), \ell_t(v))=(m-1,m)$ or $(\ell_t(u), \ell_t(v))=(m, m-1)$, we disconnect $e$ by duplicating the vertex $v$ and separating the resulting copies (one copy remaining attached to the parent, the other to the descendants). 
    We repeat this procedure for every edge of type $m-1 \rightarrow m$ or $m \rightarrow m-1$. We obtain in this way a collection of subtrees. One of these is the \textit{root component}, namely the connected component containing the root of $t$.  The other components are rooted either at label $m$ or at label $m-1$. In the first case, subtracting $m-1$ from every label in the connected component yields an element of $\mathcal E_+$. In the second case, subtracting $m$ to every label yields an element of $\mathcal E_-$. The number of positive tree excursions (resp. the number of negative excursions) obtained in this way is $e_{m-1, m}^t$ (resp. $e_{m, m-1}^t$). We denote by
    \begin{align*}
    	\{\tau_i^{m, +}\}_{i\leq e_{m-1, m}^t}=\{\tau_1^{m, +}, \dots, \tau_{e_{m-1, m}^t}^{m, +}\},
    \end{align*}
    the unordered collection of these positive tree excursions, and by
    \begin{align*}
    	\{\tau_j^{m, -}\}_{j\leq e^t_{m, m-1}}=\{\tau_1^{m, -}, \dots, \tau_{e_{m, m-1}^t}^{m, -}\},
    \end{align*}
    the unordered collection of the negative tree excursions. 
    
    \medskip
	
	In the second step (step (II) in Figure \ref{Fig1}), we construct a plane forest $\mathcal F^m_{t}$, called the \emph{excursion forest}, whose vertex set is 
	\begin{align}
		\label{AllExcursions}
	\mathfrak E= \{\tau_i^{m, +}\}_{i\leq e_{m-1, m}^t}\cup \{\tau_j^{m, -}\}_{j\leq e^t_{m, m-1}},
	\end{align} the collection of all positive and negative tree excursions. Each component detached directly from the root component $T^{[m]}$ gives rise to the root of a tree in $\mathcal F^m_{t}$. Moreover, a tree excursion $\tau\in \mathfrak E$ is the child in $\mathcal F_t^m$ of another tree excursion $\tau'\in \mathfrak E$ if and only if the root of $\tau$ corresponds to a leaf of $\tau'$ in the initial tree $t$. In this way, we obtain a forest $\mathcal F^m_{t}$ of $|\{v \in t^{[m]}  : \ell_t(v) = m\}|$ trees, whose vertices are themselves tree excursions. The forest $\mathcal F^m_{t}$ inherits a plane forest structure from the plane tree structure of $t$.
	
\subsection{A combinatorial lemma}

We now record a counting lemma, which is a consequence of \cite[Lemma~4.1]{MarkovPropertyMD} and which will be useful later.
\begin{lemma}
	\label{lemmaCounting}
	Let $r\geq 0$, $p\geq 0$ be integers and  $(m^-_{i})_{i\geq 0}$ and $(m^+_{i})_{i\geq 0}$ be two collections of nonnegative integers such that
	\begin{align}
		\label{CondiLemma}
		r+1=\sum_{i\geq 0}m_i^-=1+\sum_{i\geq 0} i m_i^+ \qquad \text{ and } 	\qquad p=\sum_{i\geq 0} m_i^+=\sum_{i\geq 0} i m_i^-.  
	\end{align}
	Let $\mathcal{N}((m^-_{i})_i, (m^+_{i})_i)$ be the number of plane trees such that for every $i\geq 0$:
	\begin{itemize} [label=\scriptsize$\bullet$]
		\item The number of vertices at even height with $i$ children is exactly $m^-_i$.
		\item The number of vertices at odd height with $i$ children is exactly $m^+_i$.
	\end{itemize}
	Then, we have
	\begin{align}
		\label{CountFormulaN}
		\mathcal{N}\left((m^-_{i})_{i}, (m^+_{i})_{i} \right) = \frac{p!r!}{\displaystyle \prod_{i\geq 0}m^-_{i}!\prod_{i\geq 0}m^+_{i}!}.
	\end{align}
\end{lemma}
\begin{proof}[Proof of Lemma \ref{lemmaCounting}]
	Let $t$ be a plane tree having exactly
	$m_i^-$ (resp. $m_i^+$) vertices at even (resp. odd) height with $i$ children, for every integer $i \geq 0$.
	Thus, by assumption \eqref{CondiLemma}, $t$ has $\sum_{i\geq 0}m_i^-=r+1$ vertices at even height and $\sum_{i\geq 0} m_i^+=p$ vertices at odd height.
	Let $u_0$ denote the root of $t$ and label the non-root vertices at even height by
	$u_1,\ldots,u_r$ and the vertices at odd height by
	$v_1,\ldots,v_p$. Let $n^+_k$ (resp. $n^-_j$) be the number of children of $u_k$ for
	$0\leq k\leq r$ (resp. of  $v_j$ for $1\leq j\leq p$). Then, for every $i\geq0$, we have
	\[ \left|\{j \in [\![1,p]\!]:n^-_j=i\} \right|=m_i^+
	\qquad \text{and} \qquad
	\left|\{k \in [\![0,r]\!]:n^+_k=i\}\right|=m_i^-. \]
	Consequently, we deduce that \[r+1=1+\sum_{i\geq0} i\,m_i^+
	=1+\sum_{j=1}^{p} n^-_j,
	\qquad \text{and} \qquad 
	p=\sum_{i\geq0} i\,m_i^-
	=\sum_{k=0}^{r} n^+_k.\]
	Viewing even-height vertices as positive and odd-height vertices as negative, the assumptions of \cite[Lemma 4.1]{MarkovPropertyMD} are satisfied with
	$n=1$, $p=r+1$ and $q=p$.  Hence, \cite[Lemma 4.1]{MarkovPropertyMD} implies that for fixed labels $u_1,\ldots,u_r$ and $v_1,\ldots,v_p$, the number of such trees is $ p!\,r!.$ \medskip
	
	\noindent It remains to remove the labels. Permuting vertices having the same number of children does not
	change the underlying plane tree. For every $i\geq0$, there are
	$m_i^-!$ ways to permute the even-height vertices with $i$ children and
	$m_i^+!$ ways to permute the odd-height vertices with $i$ children.
	Therefore each underlying plane tree is counted exactly $\prod_{i\geq0}m_i^-!\prod_{i\geq0}m_i^+!$ times among the labelled trees counted above.
	Hence, dividing by this multiplicity yields
	\[\mathcal N((m^-_{i})_i,(m^+_{i})_i)
	=
	\frac{p!\,r!}
	{\displaystyle
		\prod_{i\geq0}m_i^-!\prod_{i\geq0}m_i^+!},\]
	which completes the proof of Lemma \ref{lemmaCounting}.
\end{proof}

\subsection{Markov properties}
Let $T$ be sampled with law $\Pi_0$. Recall that we defined, for every $k\geq 1$, 
\begin{align*}
	\mathcal Y^+_{k}=(e^T_{k-1, k}, e^T_{k-1, k-1}, e^T_{k, k-1})\qquad \text{ and } \qquad 	\mathcal Y^+_{k}=(e^T_{-k+1, -k}, e^T_{-k+1, -k+1}, e^T_{-k, -k+1}).
\end{align*}
The following proposition is an adaptation of results in \cite{MarkovPropertyMD}.
\begin{proposition}\label{MP} 	Under $\Pi_0$, the processes $\mathcal Y^+$ and $\mathcal Y^-$  are time-homogeneous Markov chain under $\Pi_0$ which are independent conditionally on the values taken by $\mathcal Y^+_1$ and $\mathcal Y^-_1$. Moreover,
	\begin{enumerate}
		\item 	Let $(p, q, r)\in \mathbb Z^3_{\geq 0}$, and assume that there exists $m\geq 1$ with $\Pi_0\left( {\mathcal Y}_m^{+}=(p, q, r)\right)>0$. Let $\mathtt t_1, \cdots, \mathtt t_p$ be independent positive tree excursions with law $\Pi^+$. The event $\{\sum_{i=1}^p n_{\mathtt t_i}= r\}$ has positive probability and for every $(a, b, c)\in \mathbb Z_{\geq 0}^3$, we have
		\begin{align*}
			\Pi_0\left(\mathcal Y^+_{m+1}=(a, b, c) \mid \mathcal Y^+_m=(p, q, r)\right)= \mathbb P\left( \sum_{i=1}^p (e_{1, 2}^{\mathtt t_i}, e_{1, 1}^{\mathtt t_i}, e_{2, 1}^{\mathtt t_i})=(a, b, c) \ \Big| \ \sum_{i=1}^pn_{\mathtt t_i}= r\right).
		\end{align*}
		\item Symmetrically, let $(p, q, r)\in \mathbb Z^3_{\geq 0}$, and assume that $\Pi_0\left( {\mathcal Y}_m^{+}=(p, q, r)\right)>0$ for some $m\geq 1$. Let $\mathtt t_1, \cdots, \mathtt t_p$ be independent negative tree excursions with law $\Pi^-$. The event $\{\sum_{i=1}^p n_{\mathtt t_i}= r\}$ has positive probability and for every $(a, b, c)\in \mathbb Z_{\geq 0}^3$, we have
		\begin{align*}
				\Pi_0\left(\mathcal Y^-_{m+1}=(a, b, c) \mid \mathcal Y^-_m=(p, q, r)\right)=\mathbb P\left( \sum_{i=1}^p (e_{-1, -2}^{\mathtt t_i}, e_{-1, -1}^{\mathtt t_i}, e_{-2, -1}^{\mathtt t_i})=(a, b, c) \ \Big|\ \sum_{i=1}^pn_{\mathtt t_i}= r\right).
		\end{align*}
	\end{enumerate}	
\end{proposition}
\begin{remark} Note that the right-hand-sides of the transition kernel in point $1$ and $2$ of Proposition \ref{MP} do not depend on $q$.
\end{remark}
\begin{proof} Fix $m\geq 1$, and recall the notation $T^{[m]}$ for the tree $T$ truncated at label $m$. Fix moreover $(p, q, r)\in \mathbb Z^3_{\geq 0}$ and $m\geq 1$ with 
	\begin{align*}
		\Pi_0\left( {\mathcal Y}_m^+=(p, q, r)\right)>0.
	\end{align*}
	First, we argue conditionally on the event
	\begin{align*}
		\mathcal A_{p, r}^m:= \left\{(e^T_{m-1, m}, e^T_{m, m-1})= (p, r)\right\}.
	\end{align*}
	Recall the notation $\{\tau_j^{m, -}\}_{j\leq r}$ (resp. $\{\tau_i^{m, +}\}_{i\leq p}$) for the unordered collection of negative (resp. positive) tree excursions associated with the excursion forest decomposition of $T$ at level $m$. By \cite[Theorem 4.2]{MarkovPropertyMD}, conditionally on $\mathcal A_{p, r}^m$, the unordered collection $\{\tau_i^{m, +}\}_{i\leq p}$ is independent of both $T^{[m]}$ and $\{\tau_j^{m, -}\}_{j\leq r}$. Every edge below label $m-1$ belongs to the truncation $T^{[m]}$ or to one of the negative tree excursions. For every $(i, j)\in \Delta$ with $i, j\leq m-1$, we have
	\begin{align}
		\label{FormulaBelowLabelm}
		e^T_{i, j}= e^{T^{[m]}}_{i, j}+ \sum_{k=1}^r e^{\tau_k^{m, -}}_{i-m, j-m}.
	\end{align}
 In particular $e_{m-1, m-1}^T$ is a measurable function of  $T^{[m]}$ and $\{\tau_j^{m, -}\}_{j\leq r}$ and therefore is independent of $\{\tau_i^{m, +}\}_{i\leq p}$ conditionally on $\mathcal A_{p, r}^m$. It follows that the conditional law of $\{\tau_i^{m, +}\}_{i\leq p}$ given $\mathcal A_{p, r}^m$ is the same as the conditional law of $\{\tau_i^{m, +}\}_{i\leq p}$ given $\{ {\mathcal Y}_m^+=(p, q, r)\}$.  Moreover, conditionally on $\{{\mathcal Y}_m^+= (p,q, r)\}$, the collection $\{\tau_i^{m, +}\}_{i\leq p}$ is independent of $T^{[m]}$ and $\{\tau_j^{m, -}\}_{j\leq r}$.\\
	
 We now argue conditionally on $\{\mathcal Y_m^+=(p, q, r)\}$. By \eqref{FormulaBelowLabelm}, we see that every ${\mathcal Y}_l^+$ with $l< m$ is a measurable function of $T^{[m]}$ and $\{\tau_j^{m, -}\}_{j\leq r}$. For every $(i, j)\in \Delta$ with $i, j\geq m$, we have:
	\begin{align}
		\label{FormulaAboveLabelm}
		e^T_{i, j}= \sum_{k=1}^p e^{\tau_k^{m, +}}_{i-m+1, j-m+1}.
	\end{align}
	This yields that ${\mathcal Y}_{m+1}^+$ is a measurable function of $\{\tau_i^{m, +}\}_{i\leq p}$. It follows from the previous discussion that $ {\mathcal Y}_{m+1}^+$ is independent of $( {\mathcal Y}_l^+)_{l<m}$ conditionally on $ {\mathcal Y}_m^+$.  This exactly means that the process ${\mathcal Y}^+$ 
    is a Markov chain. \medskip
   
	To see that this chain is time-homogeneous, recall from the previous discussion that the law of $\mathcal Y_{m+1}^+$ conditionally on $\{{\mathcal Y}_m^+= (p,q, r)\}$ is the same as the law of $\mathcal Y_{m+1}^+$ conditionally on $\mathcal A_{p, r}^m$. According to \cite[ Theorem 4.2]{MarkovPropertyMD}, conditionally on $\mathcal A_{p, r}^m$, the collection $\{\tau_i^{m, +}\}_{i\leq p}$ is distributed as $p$ independent positive tree excursions with law $\Pi^+$, conditioned to have exactly $r$ leaves labelled $0$. More precisely,  if $\mathcal A_{p, r}^m$ has positive probability, then if we let $(\mathtt t_1, \cdots,\mathtt t_p)$ be a collection of $p$ independent tree excursions with law $\Pi^+$, the event
\begin{align}
	\label{Eventsumnti}
	\Big\{\sum_{i=1}^p n_{\mathtt t_i}= r\Big\},
\end{align}
 has positive probability, where we recall that $n_{\mathtt t_i}$ is the number of leaves labelled $0$ in $\mathtt t_i$. Moreover, conditionally on $\mathcal A_{p, r}^m$ the law of $\{\tau_i^{m, +}\}_{i\leq p}$ is the same as the law of $\{\mathtt t_i\}_{i\leq p}$ conditionally on the event in \eqref{Eventsumnti}. This implies that for every $a, b, c\in \mathbb Z_{\geq 0}$:
	\begin{align*}
		\Pi_0\left({\mathcal Y}_{m+1}^+=(a, b, c)  \ \Big|\  {\mathcal Y}^+_m=(p, q, r) \right)=\mathbb P\left( \sum_{i=1}^p (e_{1, 2}^{\mathtt t_i}, e_{1, 1}^{\mathtt t_i}, e_{2, 1}^{\mathtt t_i})=(a, b, c) \ \Big| \ \sum_{i=1}^pn_{\mathtt t_i}= r\right).
	\end{align*}
	 This proves point $2$ in Proposition \ref{MP}. Moreover since the right-hand-side is independent of $m$, this proves that $({\mathcal Y}_{m}^+)_{m\geq 1}$ is a \emph{time-homogeneous} Markov chain. Point $3$ in Proposition \ref{MP} follows by similar arguments and shows that $({\mathcal Y}_m^-)_{m\geq 1}$ is also a time-homogenous Markov chain.
	
	\medskip 
	
	It remains to show that $\mathcal Y^-$ and $\mathcal Y^+$ are independent conditionally on the values taken by $\mathcal Y^+_1$ and $\mathcal Y^-_1$. Using \eqref{FormulaAboveLabelm} with $m=1$, we note that the process $({\mathcal Y}_{l}^+)_{l\geq 1}$ is a measurable function of $\{\tau_i^{1, +}\}_{i\leq e_{0, 1}^T}$ and ${\mathcal Y}_1^+$. Similarly, using \eqref{FormulaBelowLabelm} with $m=1$, we get that $({\mathcal Y}_{l}^-)_{l\geq 1}$ is a measurable function of $T^{[1]}$ and $\{\tau_j^{1, -}\}_{j\leq e^T_{1, 0}}$. As a direct consequence of \cite[Theorem 4.2]{MarkovPropertyMD} and the previous discussion for $m=1$, we get that these two processes are independent conditionally on the values taken by ${\mathcal Y}_1^+$ and ${\mathcal Y}_1^-$. This completes the proof of the proposition.
\end{proof}

\subsection{The excursions below label $1$.} 
Let $\tau_0^{1, -}$ be the tree in $\mathbb T_{-1}$ obtained by shifting all labels of $T^{[1]}$ by $-1$. In particular $\tau_0^{1, -}\in \mathcal E_-$. Recall that $(\tau_i^{1, -}; 1\leq i \leq e^T_{1, 0})$, are the tree excursions below label $1$ in the forest excursion decomposition. The following proposition describe the law of 
\begin{align*}
	\left\{\tau_0^{1, -}, \cdots, \tau_{e^T_{1, 0}}^{1, -}\right\},
\end{align*}
conditionally on the values taken by $e^T_{0, 1}$ and $e^T_{1, 0}$. Its proof follows the arguments in the proof of \cite[Theorem~4.2]{MarkovPropertyMD}.
\begin{proposition}
			\label{lemmaRootTransition}
			Let $T$ be sampled with law $\Pi_0$. Fix $(p,r)\in \mathbb Z_{\geq 0}^2$ such that the event
			\begin{align}
				\label{EventCondi}
				\{e_{0, 1}^T= p, \ \ e_{1, 0}^T= r\},
			\end{align}
			has positive probability. Conditionally on the event in \eqref{EventCondi}, the collection $\{\tau_0^{1, -}, \cdots, \tau_{r}^{1, -}\}$ of negative tree excursions has the law of an unordered collection $\{\mathtt t_0, \cdots, \mathtt t_r\}$ of $r+1$ independent trees with law $\Pi^-$, conditioned on the event
			\begin{align*}
				\left\{n_{\mathtt t_0}+ \cdots+ n_{\mathtt t_r}= p\right\}.
			\end{align*}
		\end{proposition}

\begin{proof}[Proof of Proposition \ref{lemmaRootTransition}] Fix $(p,r)\in \mathbb Z_{\geq 0}^2$ such that the event in \eqref{EventCondi} has positive probability. In particular, note that this imply that $r=0$ if $p=0$, since we must have at least one edge from label $0$ to $1$ to have an edge from label $1$ to $0$. The result of Proposition \ref{lemmaRootTransition} is automatically satisfied when $p=r=0$. We now assume that $p>0$ and we argue conditionally on the event in \eqref{EventCondi}. For every $k\geq 0$ and every tree excursions $t,t_1, t_2, \cdots, t_k$ in $\mathcal E_+$ (or $\mathcal E_-$), write
	\begin{align}
		\label{nboccurences}
		C_t(t_1, \cdots, t_k)= \mathrm{Card}\{i\leq k \mid t_i=t\},
	\end{align} corresponding to the number of occurrences of $t$ in the sequence $t_1, \cdots, t_k$.  For every $t\in \mathcal E_+$, let $\mathtt c_t^{1,+} \coloneqq C_t(\tau_1^{1,+}, \cdots, \tau_p^{1,+})$ and for every $t\in \mathcal E_-$, let $\mathtt c_t^{1,-} \coloneqq C_t(\tau_1^{1,-}, \cdots, \tau_r^{1,-})$. \\
    
  Fix $t_0 \in \mathcal E_-$ and two collections $(c_t^+)_{t\in \mathcal E_+}$ and $(c_t^-)_{t\in \mathcal E_-}$ of non-negative integers with finitely many non-zero terms, such that
	\begin{align*}
		p = \sum_{t\in \mathcal E_+} c_t^+ = n_{t_0} + \sum_{t\in \mathcal E_-} c_t^- n_t,
		\qquad
		r = \sum_{t\in \mathcal E_-} c_t^- = \sum_{t\in \mathcal E_+} c_t^+ n_t.
	\end{align*}
	 Write $t_0'$ for the element of $\mathbb T_0$ obtained by shifting every label of $t_0$ by $+1$. According to Equation (18) in \cite{MarkovPropertyMD}, the probability $\Pi_0(T^{[1]}=t_0', (\mathtt c_t^{1,+})_{t\in \mathcal E_+} \!= (c_t^{+})_{t\in \mathcal E_+},  (\mathtt c_t^{1,-})_{t\in \mathcal E_-} \!= (c_t^{-})_{t\in \mathcal E_-} )$ is equal to
     \begin{align*}
		n_{t_0} r!(p-1)! \, \Pi^{[1]}_0(t_0') \, \prod_{t\in \mathcal E_-} \frac{\Pi^{-}(t)^{c_t^-}}{c_t^-!} \, \prod_{t\in \mathcal E_+} \frac{\Pi^{+}(t)^{c_t^+}}{c_t^+!}.
	\end{align*}
	Using the definition of $\tau_0^{1,-}$ and the identity $\Pi_0^{[1]}(t_0') = \Pi^{-}(t_0)$, we deduce that the probability $\Pi_0(\tau_0^{1,-}=t_0, \ (\mathtt c_t^{+})_{t\in \mathcal E_+}= (c_t^{+})_{t\in \mathcal E_+}, \ (\mathtt c_t^{-})_{t\in \mathcal E_-}= (c_t^{-})_{t\in \mathcal E_-} )$ is equal to
     \begin{equation}\label{EqualityAround0}
		n_{t_0} r!(p-1)! \, \Pi^{-}(t_0)\, \prod_{t\in \mathcal E_-} \frac{\Pi^{-}(t)^{c_t^-}}{c_t^-!} \, \prod_{t\in \mathcal E_+} \frac{\Pi^{+}(t)^{c_t^+}}{c_t^+!}.
	\end{equation}
For every $t\in \mathcal E_-$, define $\widetilde{\mathtt c}_t^{1,-} = \mathtt c_t^{1,-} + 1$ if $t = \tau_0^{1,-}$ and $\widetilde{\mathtt c}_t^{1,-} = \mathtt c_t^{1,-}$ otherwise. Equivalently,
	\[\widetilde{\mathtt c}_t^{1,-} = C_t(\tau_0^{1,-},\tau_1^{1,-}, \cdots, \tau_r^{1,-}).\] 
Similarly, set $\widetilde c_{t_0}^- = c_{t_0}^- + 1$ and $\widetilde c_t^- = c_t^-$ for $t \neq t_0$. From \eqref{EqualityAround0}, we obtain that the probability $\Pi_0(\tau_0^{1,-}=t_0, \ (\mathtt c_t^{+})_{t\in \mathcal E_+}= (c_t^{+})_{t\in \mathcal E_+}, \ (\widetilde{\mathtt c}_t^{-})_{t\in \mathcal E_-}= (\widetilde c_t^{-})_{t\in \mathcal E_-})$ is equal to
	\begin{align}
		 n_{t_0} \widetilde c_{t_0}^- r!(p-1)! \prod_{t\in \mathcal E_-} \frac{\Pi^{-}(t)^{\widetilde c_t^-}}{\widetilde c_t^-!} \prod_{t\in \mathcal E_+} \frac{\Pi^{+}(t)^{c_t^+}}{c_t^+!}.  \label{EqualityAround02}
	\end{align}
Now, fix collections $(d_t^+)_{t\in \mathcal E_+}$ and $(\widetilde d_t^-)_{t\in \mathcal E_-}$ of non-negative integers such that
\begin{align*}
		p = \sum_{t\in \mathcal E_+} d_t^+ = \sum_{t\in \mathcal E_-} \widetilde d_t^- n_t,
		\qquad
		r = \sum_{t\in \mathcal E_-} \widetilde d_t^- - 1 = \sum_{t\in \mathcal E_+} d_t^+ n_t.
	\end{align*}
	Decomposing over all events $\{\tau_0^{1,-}=t_0\}$ for $t_0\in \mathcal E_-$ and using \eqref{EqualityAround02}, we get 
\begin{align*}
		& \Pi_0\left((\mathtt c_t^{+})_{t\in \mathcal E_+}= (d_t^{+})_{t\in \mathcal E_+}, \ (\widetilde{\mathtt c}_t^{-})_{t\in \mathcal E_-}= (\widetilde d_t^{-})_{t\in \mathcal E_-} \right) \\
		&  =\sum_{t_0\in \mathcal E_-}\Pi_0\left(\tau_0^{1,-}=t_0, \ (\mathtt c_t^{+})_{t\in \mathcal E_+}= (d_t^{+})_{t\in \mathcal E_+}, \ (\widetilde{\mathtt c}_t^{-})_{t\in \mathcal E_-}= (\widetilde d_t^{-})_{t\in \mathcal E_-} \right) \\
		&  = \left(\sum_{t_0\in \mathcal E_-} n_{t_0}\widetilde d_{t_0}^- \right) r!(p-1)!  \prod_{t\in \mathcal E_-} \frac{\Pi^{-}(t)^{\widetilde d_t^-}}{\widetilde d_t^-!} \, \prod_{t\in \mathcal E_+} \frac{\Pi^{+}(t)^{d_t^+}}{d_t^+!} \\
		& = r!p! \prod_{t\in \mathcal E_-} \frac{\Pi^{-}(t)^{\widetilde d_t^-}}{\widetilde d_t^-!} \, \prod_{t\in \mathcal E_+} \frac{\Pi^{+}(t)^{d_t^+}}{d_t^+!}
	\end{align*}
where the last equality follows from the identity $\sum_{t_0\in \mathcal E_-} n_{t_0}\widetilde d_{t_0}^- =p$. Note the symmetry between the law of $(\mathtt c_t^+)_{t\in \mathcal E_+}$ and $(\widetilde{\mathtt c}_t^-)_{t\in \mathcal E_-}$. In particular, we can follow the same arguments as in the end of the proof of Theorem 4.2 of \cite{MarkovPropertyMD} to get that, conditionally on  the event in \eqref{EventCondi},
	\begin{align*}
		(\widetilde{\mathtt c}_t^-)_{t\in \mathcal E_-} \overset{(d)}{=} (C_t(\mathtt t_0, \cdots, \mathtt t_r))_{t\in \mathcal E_-},
	\end{align*}
	where $\mathtt t_0, \ldots, \mathtt t_r$ is a family of $r+1$ negative tree excursions with law $\Pi^-$ conditioned on having a total of $p$ leaves labelled $0$. This completes the proof of Proposition \ref{lemmaRootTransition}.
\end{proof}
    
	\section{Some explicit formulas}
	In this section, we derive explicit formulas 
    related to the vertical profile of our labelled tree models. We first obtain expressions for the transition kernels of the Markov chains $\mathcal Y^+ $ and $\mathcal Y^-$ in terms of convolutions of the laws $\nu_+$, $\nu_-$ and $\eta_*$ introduced in Section \ref{KernelFormSec}. 
    We also characterize the law of the vertical profile at $0$, through combinatorial and probabilistic identities involving Bienaymé--Galton--Watson trees and their {\L}ukasiewicz walks. Finally, combining these results, we establish a product formula for the probability that a tree has a prescribed vertical edge profile.
    
 \subsection{Proof of the transition kernel formulas (Theorem \ref{FormKernel})} \label{sect_kernel}

\begin{proof}[Proof of Theorem \ref{FormKernel}]  Recall from section \ref{KernelFormSec} that ${\mathbf R}$ and $\mathbf L$ denote the transition kernels of ${\mathcal Y}^+$ and $\mathcal Y^-$. 
Fix $p, q, r\geq 0$ with $\Pi_0({\mathcal Y}^+_m=(p, q, r))>0$ for some $m\geq 1$. Fix $a, b, c\geq 0$. By definition
 \begin{align*}
 	{\mathbf R}_{p, q, r}^{a, b, c}= \Pi_0\left(  {\mathcal Y}^+_{m+1}= (a, b, c) \mid {\mathcal Y}^+_m=(p,q, r)\right).
 \end{align*}
 Let $t_1, \cdots,  t_p$ be independent tree excursions with law $\Pi_+$. In particular the numbers $n_{t_i}, i\leq p$ of leaves labelled $0$ in each $t_i$ are independent with law $\nu_+$. 	If $t\in \mathcal E_+$, we write
 \begin{align*}
 	y_t^+ \coloneqq e_{1, 2}^t, \qquad y_t^{\circ} \coloneqq e_{1, 1}^t, \qquad \text{and} \qquad  y_t^-\coloneqq e_{2, 1}^t,
 \end{align*} to lighten notation.
By point $1$ in Proposition \ref{MP}, we have $
  f^{\nu_+}_p(r)=	\mathbb P(\sum_{i=1}^pn_{t_i}=r)>0$,
  and
 \begin{align}
 		{\mathbf R}_{p, q, r}^{a, b, c}&= \mathbb P\left( \sum_{i=1}^p (e_{1, 2}^{ t_i}, e_{1, 1}^{ t_i}, e_{2, 1}^{ t_i})=(a, b, c) \ \Big| \ \sum_{i=1}^pn_{ t_i}= r\right)\nonumber \\
 		&= \mathbb P\left(\sum_{i=1}^p(y_{t_i}^+, \ y_{t_i}^\circ,\  y^-_{t_i})  = (a , b, c) \ \Big|\ \sum_{i=1}^p n_{t_i}= r \right)\nonumber \\
 		&=\label{Midequation}
		 \frac{1}{f^{\nu_+}_p(r)}\mathbb P\left(\sum_{i=1}^py_{t_i}^+ = a \ \ , \ \sum_{i=1}^py_{t_i}^\circ= b  \ , \ \  \sum_{i=1}^p y^-_{t_i}= c \ , \ \  \sum_{i=1}^p n_{t_i}= r \right).
\end{align}
	If $t \in \mathcal{E}_+$ is a positive tree excursion with law $\Pi^+$, we introduce the tree $R_t$ \emph{induced by the vertices labelled $1$}. More precisely, write $\mathcal V_t=\{u\in V(t)\mid \ell_t(u)=1\}$ for the subset of vertices of $t$ with label $1$. The tree $R_t$ has vertex set $\mathcal V_t$ and for every $u, v\in \mathcal V_t$, an edge is drawn from $u$ to $v$ in $R_t$ if and only if $u$ is a strict ancestor of $v$ in $t$, and no other vertex of $\mathcal V_t$ lies strictly between $u$ and $v$ along the ancestral line of $v$. The tree $R_t$ inherits a natural plane tree structure from $t$. If $u\in \mathcal V_t$, let us moreover write:
	\begin{itemize}[label=\scriptsize$\bullet$]
		\item $\iota_t^-(u)$ for the number of children of $u$ with label $0$ in $t$,
		\item $\iota_t^\circ(u)$ for the number of children of $u$ with label $1$ in $t$,
		\item $\iota_t^+(u)$ for the number of children of $u$ with label $2$ in $t$.
	\end{itemize}
	\begin{figure}[h] 
			\centering
			\includegraphics[scale = 1]{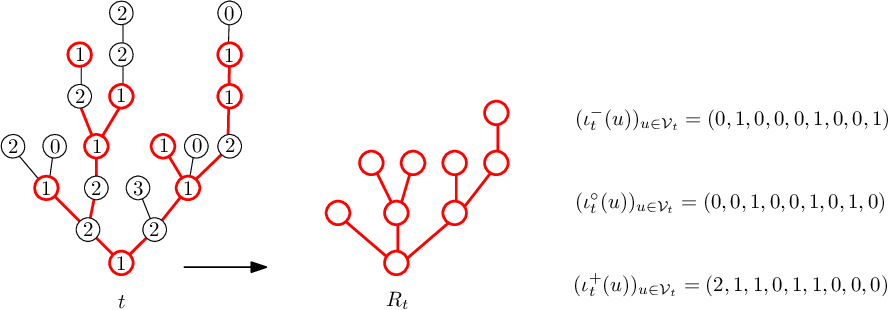}
			\caption{An example of a tree excursion $t$, its reduced tree $R_t$ and the values taken by the labelling functions, in the lexicographic enumeration of the vertices of $R_t$.}
			\label{Fig2}
		\end{figure}
		By construction, we have:
	\begin{align}
		\label{yfromFrakTch4}
		|R_t| = y^-_t+y^\circ_t, \qquad \sum_{u\in \mathcal V_t} \iota_t^-(u) = n_t,  \qquad \sum_{u\in \mathcal V_t}\iota_t^\circ(u) = y^\circ_t, \qquad \sum_{u\in \mathcal V_t}\iota_t^+(u) = y^+_t,
	\end{align}
	where we recall that $|R_t|$ is the number of edges in $R_t$.
	The second, third and fourth equalities follow directly from the definitions of $\iota_t^-$, $\iota^\circ_t$ and $\iota^+_t$. For the first equality, observe that every vertex in $\mathcal{V}_t$, except the root, is connected to its parent in $t$ either by an oriented edge from label $2$ to label $1$ (which contributes to $y^-_\tau$) or by an oriented edge from label $1$ to label $1$ (which contributes to $y^\circ_\tau$). The first equality then follows from the fact that we have $\mathrm{Card}(\mathcal V_t)= 1+|R_t|$, since $R_t$ is a tree. The tuple: 
	\begin{align*}
		\biggr(R_t,\left(\iota_t^-(u),\iota_t^\circ(u), \iota_t^+(u)\right)_{u\in \mathcal V_t}\biggr),
	\end{align*} can be thought of as a random variable taking value in the set of plane trees with a three-dimensional label function from the set of vertices into $(\mathbb Z_{\geq 0})^3$.   Write $\mathfrak P$ for the law of this random variable.\\

	\noindent Due to the inductive nature of the construction of our models of labelled Bienaymé--Galton--Watson trees, the reduced tree also has an inductive description. If we consider a vertex $u$ of $R_t$, the labels $(\iota^-_t(u), \iota_t^\circ(u), \iota_t^+(u))$ correspond to the number of children of $u$ in $t$, with label $0$, $1$, and $2$ respectively. In particular this triplet of labels has law $\eta_*$ by definition. Each of the $\iota_t^\circ(u)$ children of $u$ with label $1$ will contribute to one child of $u$ in $R_t$, and each of the $\iota_t^+(u)$ children of $u$ with label $2$ will generate a random number with law $\nu_+$ of children of $u$ in $R_t$. All in all, we deduce that a random labelled tree 
	\begin{align*}
		(R,(l^-(u), l^\circ(u), l^+(u))_{u\in V(R)}),
	\end{align*} with law $\mathfrak P$ can be constructed inductively as follow:
	\\
	
	\noindent We start with a single vertex at generation $n=0$ and we construct the labels and the plane structure of the tree inductively. Suppose that the set of vertices at generation $n\geq 0$, as well as the labels up to generation $n-1$, have already been constructed. We construct generation $n+1$ as follows; for each vertex $u$ at generation $n$, first we sample a random variable $(l^-(u), l^\circ(u), l^+(u))$ with law $\eta_*$, independently for each $v$ and also independently of the past. The triple $(l ^-(u), l^\circ(u), l^+(u))$ will be the label of $u$.
	Conditionally on these random variables, each vertex $u$ at generation $n$ gives birth to a random number $c_{R}(u)$ of children at generation $n+1$, where:
	\begin{align}
		\label{expressionChildReducedTree}
		c_{R}(u) = l^\circ(u)+ \sum_{k=1}^{l^+(u)} \mathcal N_u^k,
	\end{align} where $\mathcal N_u^k, k=1, \cdots, l^+(u)$ are independent random variables with law $\nu_+$, sampled independently of the past. When this procedure is iterated inductively, it stops almost surely and we let $R$ be the resulting genealogical tree.
	Then the triplet $(R, (l^-, l^\circ, l^+))$ is distributed according to the law $\mathfrak{P}$.
	A rigorous proof of this claim can be derived using a stopping line argument, akin to the one found in \cite{Chauvin}. Details are omitted.\\

	Using Equations \eqref{Midequation}, \eqref{yfromFrakTch4} and the previous discussion, we get
	\begin{align}
		\label{ReformulationPpqrs}
		{\mathbf R}_{p, q, r}^{a,b,c}=	\frac{1}{f^{\nu_+}_p(r)}\mathbb P\left( \sum_{i=1}^p \sum_{u \in V(R^{(i)})} (l_i^-(u), l^\circ_i(u), l^+_i(u))= (r, b, a), \ \sum_{i=1}^p |R^{(i)}|=b+c\right),
	\end{align}
	where $(R^{(i)},  (l^-_i, l^\circ_i, l^+_i))_{i\geq 1}$ are independent random labelled trees with law $\mathfrak P$. Write $(w_1, w_2 ,\cdots)$ for the enumeration of the vertices of $(R^{(i)})_{i\geq 1}$ in depth-first-order (first enumerating all vertices of $R^{(1)}$ then all vertices of $R^{(2)}$ etc.). Fix $j\geq 1$ and let $d$ be such that $w_j\in R^{(d)}$. For sake of simplicity, write $c_j$ for the number of children of $w_j$ in the tree $R^{(d)}$ and 
	\begin{align*}
		\iota^-_j \coloneqq l^-_d(w_j), \quad \quad \ \iota^\circ_j \coloneqq l^\circ_d(w_j) \quad \text{ and } \quad \iota ^+_j \coloneqq l_d^+(w_j).
	\end{align*} It follows from the inductive description of the law $\mathfrak P$ that $(\iota_j^-, \iota_j^\circ, \iota_j^+)_{j\geq 1}$ is a sequence of independent variables with law $\eta_*$.
	Write $(S_i)_{i\geq 0}$ for the \L ukasiewicz walk associated with the enumeration $(w_j)_{j\geq 1}$ of the vertices of the forest $(R^{(i)})_{i\geq 1}$. In other words, we let $S_0=0$ and for every $i\geq 1$:
	\begin{align*}
		S_i= \sum_{j=1}^i(c_j-1).
	\end{align*}
	It follows from the description of the law $\mathfrak P$ that $(S_i)_{i\geq 0}$ is a random walk with i.i.d. steps. Moreover, conditionally on $(\iota_j^-,\iota_j^\circ,\iota_j^+)$, the random variable 
	\begin{align*}
		S_j-S_{j-1}+1-\iota_j^\circ= c_j-\iota^\circ_j,
	\end{align*}
	has the law of the sum of $\iota_j^+$ independent random variables with law $\nu_+$ (this follows from Equation \eqref{expressionChildReducedTree}). In particular, conditionally on $(\iota_j^-, \iota_j^\circ, \iota_j^+)_{j\leq n}$, we have for every $n\geq 1$ that
	\begin{align*} 
		S_n+n-\sum_{j=1}^n\iota_j^\circ= \sum_{j=1}^n(S_j-S_{j-1}+1-\iota_j^\circ),
	\end{align*}
	is distributed as the sum of $\sum_{j\leq n} \iota_j^+$ independent random variables with law $\nu_+$. That is, for every $n\geq 1$, every $\ell\in \mathbb Z$ and every $r, b, a \geq 0$, we have that
	\begin{align*}
		\mathbb P\left(S_n=\ell \ \Big| \ \sum_{j=1}^n(\iota_j^-, \iota_j^\circ, \iota_j^+)=(r, b, a)\right)= f_a^{\nu_+}(\ell+n-b).
	\end{align*}
	Moreover, since $(\iota^-_j, \iota^\circ_j, \iota^+_j)_{j\leq n}$ are i.i.d. random variables with law $\eta_*$, we have
	\begin{align*}
		\mathbb P\left( \sum_{j=1}^n(\iota_j^-, \iota_j^\circ, \iota_j^+)=(r, b, a) \right)=f^{\eta_*}_n(r, b, a).
	\end{align*}
	We finally get that
	\begin{align}
		\label{ExprJoinSch4} \mathbb P\left(S_n=\ell, \quad \sum_{j=1}^n(\iota_j^-, \iota_j^\circ, \iota_j^+)=(r, b, a)\right) =f_a^{\nu_+}(\ell+n-b)f^{\eta_*}_n(r, b, a).
	\end{align}
	Fix $p\geq 1$ and write $h_{-p}=\inf\{k\geq 0 \mid S_k=-p\}$ for the first hitting time of $-p$ by the {\L ukasiewicz} walk. The random variable $h_{-p}$ is also the first time we finish exploring the tree $R^{(p)}$ in the enumeration $(w_k)_k$ (this is a well-known fact in the \L ukasiewicz coding of forests). In other words, we have
	\begin{align*}
		h_{-p}= \sum_{k=1}^p |R^{(k)}|+p,
	\end{align*}
	since the total number of vertices in $(R^{(1)}, \cdots, R^{(p)})$ is $\displaystyle \sum_{k\leq p} |R^{(k)}|+p$.
	It follows that
    \begin{align*}
		\mathbb P  \left( \sum_{j=1}^p|R^{(j)}|=b+c, \,  \sum_{j=1}^p\sum_{u\in V(R^{(j)})} (l_j^-(u), l_j^\circ(u), l_j^+(u)) = (r, b, a)\right) 
	\end{align*}
is equal to 
\begin{align*}
    \mathbb P  \left( h_{-p}=p+b+c, \, \sum_{j=1}^{p+b+c} (\iota_j^-, \iota_j^\circ, \iota_j^+)= (r, b, a)\right).
\end{align*}
By Equation \eqref{ReformulationPpqrs}, we have
	\begin{align}
		\label{ReformulationPpqrs2}
		{\mathbf R}_{p, q, r}^{a, b, c}= \frac{1}{f^{\nu_+}_p(r)}\mathbb P\left(h_{-p}=p+b+c, \ \sum_{j=1}^{p+b+c}(\iota_j^-, \iota_j^\circ, \iota_j^+)= (r, b, a) \right).
	\end{align}
	We may now apply Kemperman's formula (see e.g. \cite[Section $6.1$]{PitmanCombin}) to obtain
    	\begin{align}\label{ReformulationPpqrs3} 
		\mathbb P \! \left( h_{-p}=p+b+c, \! \sum_{j=1}^{p\!+b\!+c}(\iota_j^-, \iota_j^\circ, \iota_j^+)= (r, b, a) \! \right) \! = \frac{p}{p+b+c}\mathbb P \! \left( \! S_{p+b+c}=-p, \! \sum_{j=1}^{p+b+c}(\iota_j^-, \iota_j^\circ, \iota_j^+)= (r, b, a) \!\right).
	\end{align}
Indeed, $(c_j, \iota_j^-, \iota_j^\circ, \iota_j^+)_{j\leq p+b+c}$ are i.i.d. random variables, and the event 
	\begin{align*}
		\left\{\sum_{j=1}^{p+b+c}(\iota_j^-, \iota_j^\circ, \iota_j^+)= (r, b, a)\right\},
	\end{align*} is invariant under cyclic permutations of these quadruples. Hence, by the classical cyclic lemma, Kemperman's formula can be applied in this case.
	We apply \eqref{ExprJoinSch4} with $\ell=-p$ and $n=p+b+c$, to the right-hand-side of \eqref{ReformulationPpqrs3} and we deduce from \eqref{ReformulationPpqrs2} that
	\begin{align*}
		{\mathbf R}_{p, q, r}^{a,b,c}	= \frac{p}{p+b+c}f^{\eta_*}_{p+b+c}(r, b, a) \frac{f_{a}^{\nu_+}(c)}{f_{p}^{\nu_+}(r)}.
	\end{align*}
Symmetrically, for every $p, q, r \geq 0$ with $\Pi_0\left({\mathcal Y}_m^{-}=(p, q, r)\right)>0$ for some $m\geq 1$, we have for every $(a, b, c)\in \mathbb Z_{\geq 0}^3$
\begin{align}\label{formulaL}
	{\mathbf L}_{p,q, r}^{a, b, c}&= \Pi_0\left(\mathcal Y^-_{m+1}=(a, b, c) \mid \mathcal Y^-_m=(p, q, r)\right)\notag \\
	&=\mathbb P\left( \sum_{i=1}^p (e_{-1, -2}^{ t_i}, e_{-1, -1}^{ t_i}, e_{-2, -1}^{ t_i})=(a, b, c) \ \Big|\ \sum_{i=1}^pn_{\mathtt t_i}= r\right) \notag\\
	  & = \frac{p}{p+b+c}f^{\eta_*}_{p+b+c}(a, b, r)\frac{f^{\nu_-}_{a}(c)}{f^{\nu_-}_{p}(r)}
\end{align}
where $( t_1, \cdots,  t_p)$ are independent negative tree excursions with law $\Pi_-$. This completes the proof of Theorem~\ref{FormKernel}.
\end{proof}	
As a direct consequence of Theorem~\ref{FormKernel}, we obtain an explicit expression for the conditional law of the local vertical edge profile around label~$0$.
\begin{corollary}
			For every $(p, r)\in \mathfrak S$ such that $\Pi_0(e_{0, 1}^T=p, e_{1, 0}^T=r)>0$, and every $a, b, c\geq0$,
			\begin{align}
				\label{c}
				\Pi_0\Big((e_{0, -1}^T, e_{0, 0}^T, e_{-1, 0}^T) =( a,  b, c) \ \Big| \ (e_{0, 1}^T, e_{1, 0}^T)=(p, r)\Big)=  \frac{(r\!+\!1)f^{\eta_*}_{c+b+r+1}(a, b, p)}{c\!+\!b\!+\!r\!+\!1}\frac{f^{\nu_-}_{a}(c)}{f^{\nu_-}_{r+1}(p)}. \quad
			\end{align}
\end{corollary}

\begin{proof}
	Fix $(p, r)\in \mathfrak S$ such that $\Pi_0(e_{0, 1}^T=p, e_{1, 0}^T=r)>0$ and $(a, b, c)\in \mathbb{Z}^3_{\geq0}$. According to Proposition \ref{lemmaRootTransition}, we have that $\Pi_0\Big((e_{0, -1}^T, e_{0, 0}^T, e_{-1, 0}^T) =( a,  b, c) \ \Big| \ (e_{0, 1}^T, e_{1, 0}^T)=(p, r)\Big)$
    is equal to 
    	\begin{align*}
		\mathbb P\left(\sum_{i=1}^{r+1} (e_{-1, -2}^{t_i}, e_{-1, -1}^{t_i}, e_{-2, -1}^{t_i})=(a, b, c) \ \Big|\ \sum_{i=1}^{r+1}e_{-1, 0}^{t_i}= p\right),
	\end{align*}
	where $(t_1, \cdots, t_{r+1})$ is a collection of $r+1$ independent tree excursions with law $\Pi^-$. Using \eqref{formulaL} with the role of $(p, r)$ played by $(r+1, p)$, we immediately get \eqref{c}.
\end{proof}

\subsection{The law of the vertical edge profile at $0$ (Theorem \ref{th_profat0})}	
This section is dedicated to the proof of Theorem \ref{th_profat0}. We start with a technical lemma.
\begin{lemma} For every $a, b, c, r\in \mathbb Z_{\geq 0}$ and every $p>0$ we have:
	\begin{align}
		\label{UsefulImplic+}
f_p^{\nu_+}(r)=0 \quad  \Rightarrow	\quad f_a^{\nu_+}(c)\, f^{\eta_*}_{p+b+c}(r, b, a)=0.
	\end{align}
	And similarly,
		\begin{align}
		\label{UsefulImplic-}
		f_p^{\nu_-}(r)=0 \quad  \Rightarrow	\quad f_a^{\nu_-}(c)\, f^{\eta_*}_{p+b+c}(a, b, r)=0.
	\end{align}
\end{lemma}
\begin{proof}
	Let us prove \eqref{UsefulImplic+}. 
    Fix $a, b, c, r\in \mathbb Z_{\geq 0}$ and $p>0$ and suppose that
	\begin{align*}
	f_a^{\nu_+}(c)>0 \quad \text{and} \quad f^{\eta_*}_{p+b+c}(r,b,a)>0,
	\end{align*}
	we will show that this implies \(f_p^{\nu_+}(r)>0\). 
    
    First, \(f_a^{\nu_+}(c)>0\) implies that there exist \(a\) tree excursions \(\tilde\tau_1,\dots,\tilde\tau_a \in \mathcal{E}_+\) with \(\Pi^+(\tilde\tau_i)>0\) for every \(i\leq a\), and
	\begin{align*}
	n_{\tilde\tau_1}+\cdots+n_{\tilde\tau_a} = c.
	\end{align*}
	Then, the condition \(f^{\eta_*}_{p+b+c}(r,b,a)>0\) ensures that there exist \(p+b+c\) trees \(t_1,\dots,t_{p+b+c} \in \mathbb{T}_1\) of height at most one, such that:
	\begin{itemize} [label=\scriptsize$\bullet$]
		\item \(\Pi_1(t_i)>0\) for every \(i\leq p+b+c\),
		 \item  the total numbers of non-root vertices with labels \(0\), \(1\), and \(2\) in  \((t_i)_{1\le i\le p+b+c}\), are \(r\), \(b\), and \(a\) respectively.
	\end{itemize}
	Let \(\tau_1,\dots,\tau_a \in \mathbb{T}_2\) be the labelled trees obtained by shifting all labels in \(\tilde\tau_1,\dots,\tilde\tau_a\) by \(+1\). 
We now construct a new forest in two steps.	First, graft to each leaf of label \(2\) in the forest \((t_i)_{i\leq p+b+c}\) one of the trees \(\tau_1,\dots,\tau_a\) by identifying its root with that leaf. This produces a forest \((t'_i)_{1\le i\le p+b+c}\) of $p+b+c$ trees in $\mathcal E_+$. By construction each \(t'_i\) has positive \(\Pi^+\)-probability, and the total number of non-root leaves labelled $0$ and $1$ in this forest are \(r\) and \(b+c\) respectively.\\

Secondly, we iteratively graft trees on leaves labelled \(1\): at each step, pick a leaf of label \(1\) in one of the trees $t_i'$ and identify it with the root of another tree in the forest. Each such operation reduces the number of trees by one and the number of leaves of label $1$ decreases by one. Since the forest initially has \(p+b+c\) trees and at least \(b+c\) leaves labelled \(1\), after exactly \(b+c\) such graftings we obtain a forest \((t''_1,\dots,t''_p)\) of \(p\) trees. Throughout this procedure, the number of leaves labelled \(0\) is unchanged, hence
	\[
	n_{t''_1}+\cdots+n_{t''_p} = r.
	\]
Moreover, each tree \(t''_i\) has positive \(\Pi^+\)-probability. Thus we have constructed \(p\) positive tree excursions whose total number of leaves is \(r\), which shows that \(f_p^{\nu_+}(r)>0\). This completes the proof of \eqref{UsefulImplic+}. The second statement \eqref{UsefulImplic-} follows by symmetry.
	\end{proof}
	
	We now prove Theorem \ref{th_profat0}, giving the laws of $(e^T_{0, 1}, e^T_{1, 0})$ and $(\mathcal Y_1^+, \mathcal Y_1^-)$.

\begin{proof}[Proof of Theorem \ref{th_profat0}] We first proceed to prove \eqref{LawiniSemi}.	Let $T$ be a labelled tree sampled according to $\Pi_0$. Write $\widetilde {\mathcal F}_{T}^1$ for the plane forest we obtain from the excursion forest ${\mathcal F}_T^1$ at level $1$ by forgetting the tree excursions associated to each vertex. We write $\mathcal T$ for the plane tree we obtained by grafting the trees of $\widetilde {\mathcal F}_T^1$ to the same root, in left-to-right order (see Figure \ref{FigRe}).\medskip
	
	\noindent
    By \cite[Lemma 3.1]{MarkovPropertyMD}, the number $N_1$ of trees in the forest $\widetilde{\mathcal F}_T^1$ has distribution $\nu_-$. Moreover, conditionally on $N_1$, we get from \cite[Proposition 3.3]{MarkovPropertyMD} that the forest $\widetilde{\mathcal F}_T^1$ has the law of a forest of $N_1$ alternating two-type Bienaymé--Galton--Watson tree with offspring distributions $(\nu_+, \nu_-)$. Consequently, the tree $\mathcal T$ has the law of an alternating two-type Bienaymé--Galton--Watson tree with offspring distributions $(\nu_-, \nu_+)$.
    That is, for every plane tree $t$, we have:
	\begin{align}
		\label{bicondiF}
		\Pi_0(\mathcal T= t)= \prod_{v\in \mathcal H_+(t)}\nu_+(c_t(v)) \prod_{u\in\mathcal H_-(t)}\nu_-(c_t(u)),
	\end{align}
	where for every plane tree $t$ we write $\mathcal H_+(t)$ (resp. $\mathcal H_-(t)$) for the set of vertices at odd height (resp. even height) in $t$.  Moreover, note that we have:
	\begin{align}
		\label{IdentifReducedtreeI}
		e^T_{0, 1}= \mathrm{Card}(\mathcal H_+(\mathcal T)), \qquad \qquad e^T_{1, 0}= \mathrm{Card}(\mathcal H_-(\mathcal T))-1.
	\end{align}	
    \begin{figure}[h]
		\centering
		\begin{minipage}[b]{0.3 \textwidth}
			\centering\includegraphics[width=\textwidth]{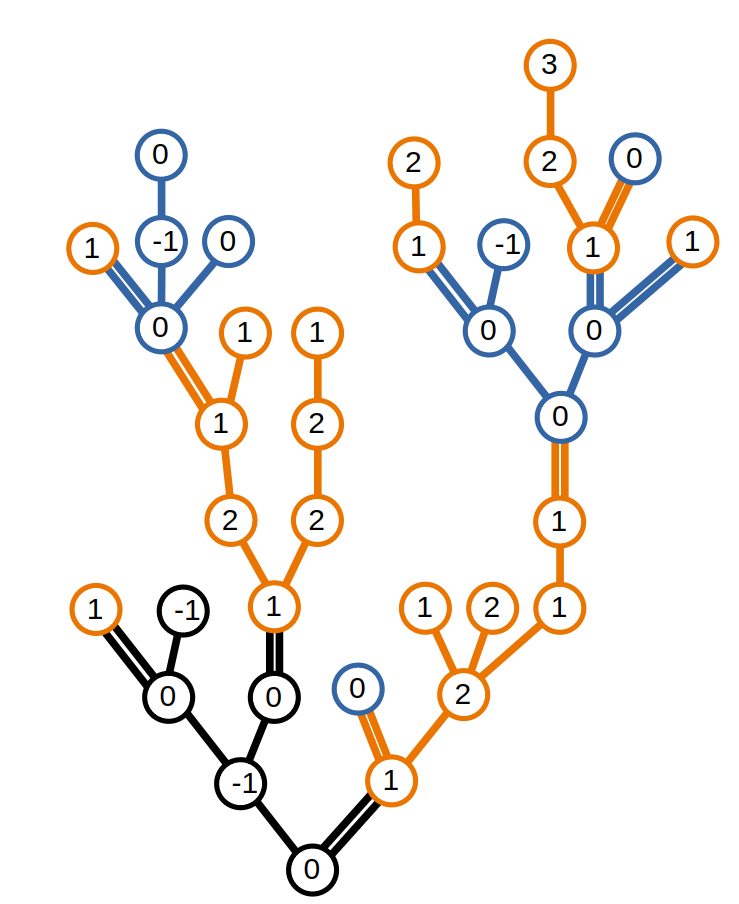}
		\end{minipage} \vrule
		\begin{minipage}[b]{0.5\textwidth}
			\centering\includegraphics[width=\textwidth]{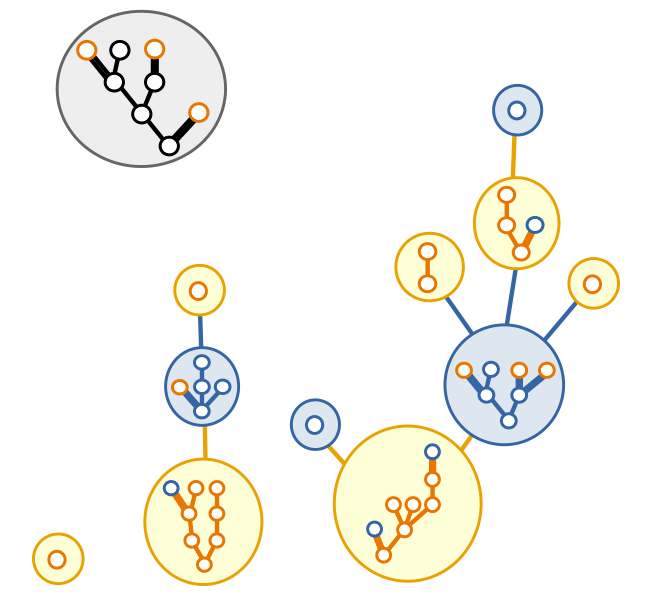}
		\end{minipage}\vrule
		\begin{minipage}[b]{0.18\textwidth}
			\centering\includegraphics[width=\textwidth]{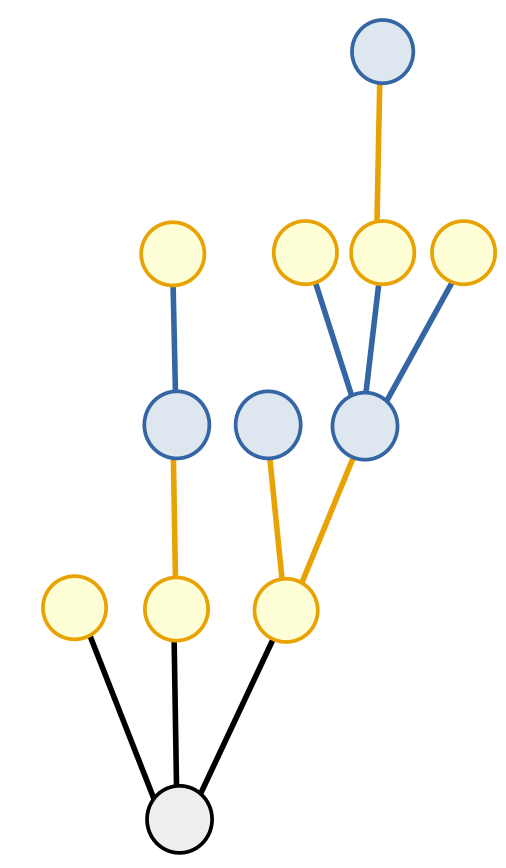}
		\end{minipage}
		\caption{
			\label{FigRe} An example of a tree $T$ (left), its forest excursion decomposition $\mathcal F^1_T$ at level $1$ and root excursion $T^{[1]}$ in grey (middle), and the reduced tree $\mathcal T$ (right). Note that the blue vertices, corresponding to the negative tree excursions below label $1$, are the non-root vertices at even height in $\mathcal T$, hence $\text{Card}(\mathcal H_-(\mathcal T))= e_{1, 0}^T+1$. Similarly the orange vertices, corresponding to the positive tree excursions above label $1$, are the vertices at odd height in $\mathcal T$, hence $\text{Card}(\mathcal H_+(\mathcal T))= e_{0, 1}^T$.}
	\end{figure}
	
	Fix $p, r\in \mathbb Z_{\geq 0}$, using \eqref{IdentifReducedtreeI} we have
	\begin{align}
		\label{TranslReducedTree} \Pi_0\left(e_{0,1}^T= p, \ e_{1,0}^T= r\right)= \Pi_0\Big(\mathrm{Card}(\mathcal H_-(\mathcal T)) = r+1, \ \mathrm{Card}(\mathcal H_+(\mathcal T))=p\Big).
	\end{align}
	Consider two collections $(m^-_{i})_{i\geq 0}$ and $(m^+_{i})_{i\geq 0}$ of nonnegative integers with finitely non-zero terms.
 Recall from Lemma \ref{lemmaCounting} the notation $\mathcal{N}((m^-_{i})_i, (m^+_{i})_i)$ for the number of plane trees $t$ such that for every $i\geq 0$,
	\begin{align*}
		\mathrm{Card}\{v\in \mathcal H_+(t) \mid c_t(v)= i\}= m_i^+, \\ \mathrm{Card}\{v\in \mathcal H_-(t) \mid c_t(v)= i\}= m_i^-.
	\end{align*}
	 Note that for such a tree $t$, we have $\mathrm{Card}(\mathcal H_+(t))= p$ and $\mathrm{Card}(\mathcal H_-(t))=r+1$ if and only if:
	\begin{align}
		\label{Condich5}
		r+1=\sum_{i\geq 0}m_i^-=1+\sum_{i\geq 0} i m_i^+ \qquad \text{ and } 	\qquad p=\sum_{i\geq 0} m_i^+=\sum_{i\geq 0} i m_i^-.  
	\end{align}
	Summing over all possible choices of $t$ and using \eqref{bicondiF} and \eqref{TranslReducedTree} we get that:
	\begin{align*}
		\Pi_0\left(e_{0,1}^T= p, \ e_{1,0}^T= r\right)= \sum_{(m^-_{i})_{i},(m^+_{i})_{i} } \mathcal{N}\left((m^-_{i})_{i}, (m^+_{i})_{i} \right) \prod_{i \geq 0} \nu_-(i)^{m^-_{i}}\nu_+(i)^{ m^+_{i}},
	\end{align*}
	where the sum is taken on all the collections $(m^-_{i})_{i\geq 0}$ and $(m^+_{i})_{i\geq 0}$ of nonnegative integers satisfying \eqref{Condich5}. Using formula \eqref{CountFormulaN} in Lemma \ref{lemmaCounting}, this yields that 
    \begin{align*}
		\Pi_0(e_{0,1}^T= p, \ e_{1,0}^T= r)&=\sum_{\substack{(m^-_{i})_{i},(m^+_{i})_{i}\\ \text{satisfying \eqref{Condich5}}}} \frac{p!r!}{\displaystyle \prod_{i\geq 0}m^-_{i}!\prod_{i\geq 0}m^+_{i}!}, \prod_{i \geq 0} \nu_-(i)^{m^-_{i}}\nu_+(i)^{ m^+_{i}}\\
		&=\frac{1}{r+1} \Biggl( \sum_{\substack{(m^+_{i})_{i} \text{ such that:} \\ \sum_i i m_i^+=r, \\ \sum_i m_i^+=p.}} \frac{ p! \underset{i\geq 0}{\prod}\nu_+(i)^{m^+_{i}}}{\underset{i\geq 0}{\prod}m_{i}^+!} \Biggr) \Biggl( \sum_{\substack{(m^-_{i})_{i} \text{ such that:}\\ \sum_i i m_i^-= p,\\ \sum_i m_i^-=r+1.}} \frac{(r+1)!\underset{i\geq 0}{\prod}{\nu}_-(i)^{m^-_{i}}}{\underset{i\geq 0}{\prod} m_{i}^-!} \Biggr).
	\end{align*}
	We now identify:
	\begin{align*}
		\sum_{\substack{(m^+_{i})_{i} \text{ such that }\\ 	r=\sum_i i\cdot m_i^+ , \ \ p=\sum_i m_i^+}} \frac{p!}{\underset{i\geq 0}{\prod}  m_{i}^+!} \prod_{i \geq 0}{\nu}_+(i)^{m^+_{i}}= [z^r]\left(\sum_{i\geq 0} \nu_+(i)z^i\right)^p= f_p^{\nu_+}(r),
	\end{align*}
	and:
	\begin{align*}
		\sum_{\substack{(m^-_{i})_{i} \text{ such that }\\ 	p=\sum_i i\cdot m_i^- , \ \ r+1=\sum_i m_i^-}} \frac{(r+1)!}{\underset{i\geq 0}{\prod}  m_{i}^-!} \prod_{i \geq 0}{\nu}_-(i)^{m^-_{i}}= [z^p]\left(\sum_{i\geq 0} \nu_-(i)z^i\right)^{r+1}= f_{r+1}^{\nu_-}(p),
	\end{align*}
	from which we finally get:
	\begin{align*}
		\Pi_0\left(e_{0,1}^T= p, \ e_{1,0}^T= r\right)= \frac{1}{r+1}f^{{\nu}_-}_{r+1}(p)f^{{\nu}_+}_{p}(r).
	\end{align*}
	This completes the proof of \eqref{LawiniSemi}.\medskip

	We now proceed to prove \eqref{LawinitialCondi}, first note that clearly we must have $b=q$ for the left-hand-side of \eqref{LawinitialCondi} to be non-zero (since ${\mathcal Y}_1^+=(p, q, r)$ implies $e^T_{0, 0}=q$ and  ${\mathcal Y}_1^{-}= (a, b, c)$ implies $e^T_{0, 0}=b$). Next, fix $p, r, a,b, c\in \mathbb Z_{\geq 0}$, we want to show that:
	\begin{align}
		\label{LawinitialCondiReform}
		\Pi_0\left((e^T_{0, -1}, e^T_{0, 0},e^T_{-1, 0}, e^T_{0, 1}, e^T_{1, 0})= (a, b, c, p, r)\right)= \frac{f^{\nu_-}_{a}(c)f^{\nu_+}_{p}(r)}{c+b+r+1}f^{\eta_*}_{c+b+r+1}(a,b,p).
	\end{align}
	We distinguish two cases
	\begin{itemize}[label=\scriptsize$\bullet$]
		\item If $\Pi_0(e^T_{0, 1}=p, e^T_{1, 0}= r)=0$, then in particular
		\begin{align*}
			\Pi_0\left((e^T_{0, -1}, e^T_{0, 0},e^T_{-1, 0}, e^T_{0, 1}, e^T_{1, 0})= (a, b, c, p, r)\right)=0.
			\end{align*}
		Moreover, using \eqref{LawiniSemi} we have $f^{{\nu}_-}_{r+1}(p)\, f^{{\nu}_+}_{p}(r)=0$. According to \eqref{UsefulImplic-} with $(r+1, p)$ playing the role of $(p, r)$, this yields
		\begin{align*}
			f^{\nu_-}_{a}(c)\, f^{\eta_*}_{c+b+r+1}(a,b,p)\, f^{{\nu}_+}_{p}(r)=0,
		\end{align*} 
		and in particular the right hand side of \eqref{LawinitialCondiReform} is zero. This completes the proof of \eqref{LawinitialCondiReform} in this case (both sides are zero).
		\item If $\Pi_0(e^T_{0, 1}=p, e^T_{1, 0}= r)>0$, we get using \eqref{c} and \eqref{LawiniSemi} that
	\begin{align*}
	\Pi_0&\left((e^T_{0, -1}, e^T_{0, 0},e^T_{-1, 0}, e^T_{0, 1}, e^T_{1, 0})= (a, b, c, p, r)\right)\\
	 &\qquad= \Pi_0\left((e^T_{0, -1}, e^T_{0, 0},e^T_{-1, 0})= (a, b, c) \mid e^T_{0, 1}=p, e^T_{1, 0}= r\right)  \Pi_0(e^T_{0, 1}=p, e^T_{1, 0}= r)\\
	 	&\qquad = \frac{(r+1)f^{\eta_*}_{c+b+r+1}(a, b, p)}{c+b+r+1}\frac{f^{\nu_-}_{a}(c)}{f^{\nu_-}_{r+1}(p)}\, \frac{1}{r+1}f^{{\nu}_-}_{r+1}(p)f^{{\nu}_+}_{p}(r)\\
	 	&\qquad = \frac{f^{\nu_-}_{a}(c)f^{\nu_+}_{p}(r)}{c+b+r+1}f^{\eta_*}_{c+b+r+1}(a,b,p).
	\end{align*}
	This shows \eqref{LawinitialCondiReform} in this case as well.
	\end{itemize}
	This completes the proof of Theorem \ref{th_profat0}.
\end{proof}

\subsection{Path probabilities}	
In this section, we use the preceding formulas to prove the product-form formula for the probability to have a prescribed vertical edge profile, as given in Proposition \ref{Path_probabilities}.
\begin{proof}[Proof of Proposition \ref{Path_probabilities}]
	For $k\geq 1$, set
	\begin{align*}
	\mathbf x_k^+ \coloneqq \big(x_{k-1,k}, x_{k-1,k-1}, x_{k,k-1}\big),
	\qquad
	\mathbf x_k^- \coloneqq \big(x_{-k+1,-k}, x_{-k+1,-k+1}, x_{-k,-k+1}\big).
	\end{align*}
	Recalling the notation \eqref{DefX+} and \eqref{DefX-}, we have
	\begin{align}
		\label{FirstExpression}
		\Pi_0\Big( e_{i,j}^T= x_{i,j}, \forall (i,j)\in \Delta\Big)
		= \Pi_0\Big({\mathcal Y}_k^+=\mathbf x_k^+, {\mathcal Y}_k^-=\mathbf x_k^-, \forall k\geq 1\Big).
	\end{align}
	 For every $m\geq 1$, we set
	\begin{align*}
		\mathcal B_m \coloneqq \Pi_0\Big(\mathcal Y_k^+=\mathbf x_k^+, \mathcal Y_k^-=\mathbf x_k^-, \forall k\leq m\Big)
	\end{align*}
	and for $m\geq 0$, we also set
	\begin{align*}
		\mathcal A_m \coloneqq
		\prod_{k=M_*\vee(-m)}^{M^*\wedge m}
		\frac{f^{\eta_*}_{n_k}(x_{k,k-1},x_{k,k},x_{k,k+1})}{n_k}
		\prod_{k=1}^{M^*\wedge m} x_{k-1,k}
		\prod_{k=M_*\vee(-m)}^{-1} x_{k+1,k}.
	\end{align*}
	This corresponds to the contribution of levels between $M_*\vee(-m)$ and $M^*\wedge m$ in the product \eqref{FormulPaths}, where we write $i\wedge j=\min(i, j)$ and $i\vee j=\max(i,j)$.
	Let us show by induction that for every $m\geq 1$,
    \begin{align}
\label{Induction}
\mathcal B_m
=
f^{\nu_+}_{x_{m-1,m}}(x_{m,m-1})
\,
f^{\nu_-}_{x_{-m+1,-m}}(x_{-m,-m+1})
\,
\mathcal A_{m-1}.
\tag{$H_m$}
\end{align}

	\noindent\emph{Case $m=1$.}
	Using Equation \eqref{LawinitialCondi}, 
	\begin{align*}
		\Pi_0\big({\mathcal Y}_1^+=\mathbf x_1^+, {\mathcal Y}_1^-=\mathbf x_1^-\big)
		&= f^{\nu_+}_{x_{0,1}}(x_{1,0})\, f^{\nu_-}_{x_{0,-1}}(x_{-1,0})\,
		\frac{f^{\eta_*}_{n_0}(x_{0,-1},x_{0,0},x_{0,1})}{n_0}= f^{\nu_+}_{x_{0,1}}(x_{1,0})\, f^{\nu_-}_{x_{0,-1}}(x_{-1,0})\,\mathcal A_0
	\end{align*}
	which establish \eqref{Induction} for $m=1$.
	
	\smallskip
	\noindent\emph{Induction step.}
	Assume \eqref{Induction} holds for some $m\geq 1$. We now prove that $(H_{m+1})$ holds by distinguishing between two cases.

	\smallskip Assume first that $\mathcal B_m>0$. Then we apply the Markov property of Proposition \ref{MP} to write:
	\begin{align}
		\label{InductionExpress} \mathcal B_{m+1}= \mathcal B_m \, {\mathbf R}^{\mathbf x_{m+1}^+}_{\mathbf x_{m}^+} \,{\mathbf L}^{\mathbf x_{m+1}^-}_{\mathbf x_{m}^-}.
	\end{align}
    Moreover, using Theorem \ref{FormKernel}
    \begin{align*}
{\mathbf R}^{\mathbf x^+_{m+1}}_{\mathbf x^+_{m}}
=
\begin{cases}
\displaystyle
x_{m-1,m}\,
\frac{
f^{\eta_*}_{n_m}
(x_{m,m-1},x_{m,m},x_{m,m+1})
}{
n_m
}
\frac{
f^{\nu_+}_{x_{m,m+1}}(x_{m+1,m})
}{
f^{\nu_+}_{x_{m-1,m}}(x_{m,m-1})
}
& \text{if } x_{m-1,m}>0,
\\[1.2em]
\mathds 1_{\mathbf x_{m+1}^+=(0,0,0)}
& \text{if } x_{m-1,m}=x_{m,m-1}=0,
\end{cases}
\end{align*}
and a similar expression holds for ${\mathbf L}_{\mathbf x_m^-}^{\mathbf x_{m+1}^-}$. Substituting these expressions into \eqref{InductionExpress}, the factors
	$f^{\nu_+}_{x_{m-1,m}}(x_{m,m-1})$ and 
	$f^{\nu_-}_{x_{-m+1,-m}}(x_{-m,-m+1})$ cancel with those in \eqref{Induction}, 
	and the remaining terms produce exactly the additional factors corresponding to levels $m$ and $-m$, yielding $(H_{m+1})$. \\
	
	  Assume now that $\mathcal B_m=0$. This implies $\mathcal B_{m+1}=0$. Using \eqref{Induction}, we distinguish three cases:
	\begin{itemize}[label=\scriptsize$\bullet$]
		\item If $\mathcal A_{m-1}=0$, then we must have $\mathcal A_m=0$ (since $\mathcal A_{m-1}$ is a factor of $\mathcal A_m$). In particular the right-hand-side of $(H_{m+1})$ is zero.
		\item  $f^{\nu_+}_{x_{m-1,m}}(x_{m,m-1})=0$, in particular $m\leq M^*$ (if $m>M^*$, we have $x_{m-1, m}=0$ and thus $x_{m, m-1}=0$ which implies 
		$f^{\nu_+}_{x_{m-1,m}}(x_{m,m-1})=1$). According to \eqref{UsefulImplic+}, we obtain
		\begin{align*}
		f^{\eta_*}_{n_m}(x_{m,m-1},x_{m,m},x_{m,m+1})\,
		f^{\nu_+}_{x_{m,m+1}}(x_{m+1,m})=0,
		\end{align*}
		and this implies that $f^{\nu_+}_{x_{m,m+1}}(x_{m+1,m})\, \mathcal A_m=0$. In particular the right-hand-side of $(H_{m+1})$ is zero.
		\item $f^{\nu_-}_{x_{-m+1,-m}}(x_{-m,-m+1})=0$, then by a symmetric argument using \eqref{UsefulImplic-} we again get that the right-hand-side of $(H_{m+1})$ is zero.
	\end{itemize} 
	 In all cases we see that $(H_{m+1})$ is satisfied (both sides are $0$). By induction, this shows that \eqref{Induction} is satisfied for every $m\geq 1$. Applying \eqref{Induction} with $m> M^*\vee(-M_*)$, we have 
	\begin{align*}
		f^{\nu_+}_{x_{m-1,m}}(x_{m,m-1})=f^{\nu_+}_0(0)=1\quad \text{and} \quad f^{\nu_-}_{x_{-m+1,-m}}(x_{-m,-m+1})=f^{\nu_-}_0(0)=1.
	\end{align*} 
 Hence $\mathcal B_m=\mathcal A_{m-1}$, which proves \eqref{FormulPaths} and completes the proof of Proposition \ref{Path_probabilities}.
\end{proof}
  
\section{Probabilistic counting of trees}
We now prove the enumeration formulas of Section \ref{ProbaCounting}. Recall the notation of Section \ref{ProbaCountingSembedded}.

\begin{proof}[Proof of Theorem \ref{theorem_embeded}  and Corollary \ref{corollaryPlanedAry}] Let $\mathcal S\subset\{-1, 0, 1\}$. Consider the model of random labelled trees associated with a critical or subcritical offspring distribution $\xi$ and uniform label increments in $\mathcal S$. That is, for every $l\geq 0$ we let:
	\begin{align*}
		\eta^{(l)}=\left(\frac{1}{|\mathcal S|}\sum_{s\in \mathcal S}\delta_s\right)^{\otimes l}.
	\end{align*}
	It follows that for every $y_{-1}, y_0, y_1, n\in \mathbb Z_{\geq 0}$ we have:
	\begin{align*}
		f_n^{\eta_*}(y_{-1}, y_0, y_1)= f^\xi_n(y_{-1}+y_0+y_1)|\mathcal S|^{-y_{-1}-y_0-y_1} \frac{(y_{-1}+y_0+y_1)!}{y_{-1}!y_0!y_1!}\1_{y_i=0,  \forall i\notin \mathcal S},
	\end{align*}
    where for every $k, n\geq 0$, $f_n^\xi(k)$ is the probability that $n$ independent random variables with law $\xi$ sum to $k$. Using Proposition \ref{Path_probabilities}, this yields:
	\begin{align*}
		&\Pi_0\left(e_{i, j}^T= x_{i, j}, \ \forall (i, j)\in \Delta_{\mathcal S}\right) = \prod_{k=1}^{M^*}x_{k-1, k} \prod_{k=M_*}^{-1} x_{k+1, k} \prod_{k=M_*}^{M^*} \frac{f^{\xi}_{n_k}(\kappa_k)|\mathcal S|^{-\kappa_k}\kappa_k!}{n_k \,x_{k, k-1}! \,x_{k,k}! \,x_{k, k+1}!}\1_{x_{k, k+s}=0, \ \forall s\notin \mathcal S},
	\end{align*} note that $\1_{x_{k, k+s}=0, \ \forall s\notin \mathcal S}$ is automatically equal to $1$ since we defined by convention $x_{i, j}=0$ whenever $(i, j)\notin\Delta_{\mathcal S}$. It follows that:
	\begin{align}
		\label{Expre_Embeddedtree}
		\Pi_0\left(e_{i, j}^T= x_{i, j} , \ \forall (i, j)\in \Delta_{\mathcal S}\right)=\frac{\displaystyle\prod_{k=1}^{M^*}x_{k-1, k} \prod_{k=M_*}^{-1} x_{k+1, k}}{\displaystyle |\mathcal S|^{N}\prod_{(i, j)\in \Delta_{\mathcal S}}x_{i, j}!} \, \prod_{k=M_*}^{M^*} \frac{\kappa_k!}{n_k}f^{\xi}_{n_k}(\kappa_k),
	\end{align}
	where we wrote $N=\sum_{k}\kappa_k$ for the total number of edges. Let us deal separately with the two cases considered in the theorem.
	\begin{enumerate}
		\item ($\mathcal A$-ary trees) Fix $\mathcal A\subset \mathbb Z_{\geq 0}$ with $0\in \mathcal A$ and $\mathcal A \neq \{0\}$. Consider the offspring distribution $\xi(k)=c_{\mathcal A}^{-1}2^{-k}\1_{k\in \mathcal A}$, where $
		c_{\mathcal A}=\sum_{k\in \mathcal A}2^{-k}$ (and again uniform displacements in $\mathcal S$).
		This offspring distribution is subcritical and we have, for every $\mathcal S$-embedded $\mathcal A$-ary tree $t$ with vertical edge profile $\underline x$,
		\begin{align*}
			\Pi_0(t)=\frac{1}{|\mathcal S|^{N}c_{\mathcal A}^{N+1}2^{N}},
		\end{align*} 
		where we wrote $N=\sum_{i, j} x_{i, j}$ for the total number of edges. In particular in this case:
		\begin{align}
			\label{CardProbaAary}
			\mathrm{Card}\left(\mathcal T_{\mathcal A\text{-ary}}^{\mathcal S}(\underline x)\right)=|\mathcal S|^{N} \,c_{\mathcal A}^{N+1} \, 2^{N}	\, \Pi_0\!\left(e_{i, j}^T= x_{i, j}, \ \forall (i, j)\in \Delta_{\mathcal S}\right).
		\end{align}
		Moreover, we easily check that $f^{\xi}_n(k)=c_{\mathcal A}^{-n}2^{-k}S_{\mathcal A}(n, k)$ in this case. Substituting this expression in \eqref{Expre_Embeddedtree} and using \eqref{CardProbaAary} yields the formula \eqref{EnumAary}. Formulas \eqref{EnumPlane} and \eqref{Enumdary} in Corollary \ref{corollaryPlanedAry} follow by replacing $S_{\mathcal A}(n ,k)$ by the explicit forms \eqref{ExplicitSAdary} and \eqref{ExplicitSplane} when $\mathcal A=\mathbb Z_{\geq 0}$ or $\mathcal A=\{0, d\}$.
		\item (Cayley trees) We use a classical link between Cayley trees and Bienaymé--Galton--Watson trees (see e.g. \cite[Example 10.2]{Janson}). Let $t$ be BGW tree associated with a Poisson offspring distribution of parameter $1$, i.e.
		\begin{align*}
			\xi(k)=\frac{1}{k!}e^{-1},
		\end{align*}
		conditioned to have $n$ vertices. If we label the vertices of $t$ by the number $1$ to $n$ uniformly at random (among the $n!$ choices) and forget the plane structure of $t$, we get a uniform rooted Cayley tree with $n$ vertices (we keep the same root as $t$). The probability for a Poisson$(1)$-BGW tree $t$ to have $n$ vertices is
		\begin{align*}
			\frac{n^{n-1}}{n!}e^{-n}.
		\end{align*}
		Hence, the probability that $t$ has $n$ vertices and the preceding procedure yields a given rooted Cayley tree with $n$ vertices is $\frac 1 {n!}e^{-n}$. So, considering again uniform displacements in $\mathcal S$, we get
		\begin{align}
			\label{CardProbCayley}
			\Pi_0\left(e_{i, j}^T= x_{i, j} , \ \forall (i, j)\in \Delta_{\mathcal S}\right)	= \frac{1}{|\mathcal S|^{N}}\frac{1}{(N+1)!} e^{-N-1} \, \mathrm{Card}\left(\mathcal T_{\text{Cayley}}^{\mathcal S}(\underline x)\right).
		\end{align}
		Moreover, we easily check that $f^{\xi}_n(k)=\frac{n^k}{k!}e^{-n}$ in this case, so substituting this expression in \eqref{Expre_Embeddedtree} and using \eqref{CardProbCayley} yields formula \eqref{EnumCayley}.
	\end{enumerate}
\end{proof}
We now prove the analogous results for canonically embedded models of trees. Recall the notation of Section \ref{CanonicallyEmbeddedSec}.
\begin{proof}[Proof of Theorem \ref{theorem_Sembedded_models}]
	We again make use of \eqref{FormulPaths} for well-chosen offspring and displacement distributions. Fix a subset $\mathcal S\subset\{-1, 0, 1\}$, write $\sigma=|\mathcal S|$ for its cardinal and assume that $\sigma \geq 2$. Let  $(s_i)_{i\leq \sigma}$ be the increasing enumeration of the elements of $\mathcal S$. To simplify notation, we write 
	\begin{align*}
		\alpha \coloneqq 1-\sigma^{-1} \qquad &\text{and} \qquad \beta \coloneqq \sigma^{-1}.
	\end{align*}
    We also wite $N=\sum_{i, j} x_{i, j}$.
	\begin{enumerate}
		\item (Complete trees) Consider the offspring distribution $\xi= \alpha\delta_0+ \beta\delta_{\sigma}$ and the displacement laws characterized by $\eta^{(|\mathcal S|)}=\delta_{(s_i)_{i\leq \sigma}}$. 
        For every canonically $\mathcal S$-embedded complete tree $t$, we have:
		\begin{align*}
			\Pi_0(t)= \beta^{\frac{N}{\sigma}}\alpha^{N+1-\frac{N}{\sigma}}.
		\end{align*}
        Note that for every $t\in \mathcal T^{\mathcal S}_{\text{comp.}}(\underline x)$, we have $|t|=N$.
		Summing over all elements of $\mathcal T^{\mathcal S}_{\text{comp.}}(\underline x)$, this yields
		\begin{align*}
			\Pi_0\left(e_{i, j}^T= x_{i, j}, \ \forall (i, j)\in \Delta_{\mathcal S}\right)= \beta^{\frac{N}{\sigma}}\alpha^{N+1-\frac{N}{\sigma}}\mathrm{Card}\left(\mathcal T^{\mathcal S}_{\text{comp.}}(\underline x)\right).
		\end{align*}
		Using \eqref{FormulPaths} it follows that:
		\begin{equation} \label{FormulPathsCompleteTree}
		\mathrm{Card}\left(\mathcal T^{\mathcal S}_{\text{comp.}}(\underline x)\right) = \beta^{-\frac{N}{\sigma}}\alpha^{-N-1+\frac{N}{\sigma}}\prod_{k=M_*}^{M^*} \left[\frac{1}{n_k}f^{\eta_*}_{n_k}(x_{k, k-1}, x_{k, k}, x_{k,k+1})\right] \prod_{k=1}^{M^*}x_{k-1, k} \prod_{k= M_*}^{-1} x_{k+1,k},
			\end{equation} 
		where we let by convention $x_{k, k+s}=0$ for every $s\notin \mathcal S$.
		Now, if $y_{-1}, y_0, y_1, n\in \mathbb Z_{\geq 0}$, note that $f^{\eta_*}_n(y_{-1}, y_0, y_1)=0$ if $y_i>0$ for some $i\notin \mathcal S$ or if $y_i\neq y_j$ for some $i, j\in \mathcal S$. Otherwise, it is immediate to check that if we let $y$ be the common value of $y_i, i\in \mathcal S$, we have:
		\begin{align*}
			f^{\eta_*}_n(y_{-1}, y_0,y_1)= \alpha^{n-y}\beta^{y}{n\choose y},
		\end{align*}
		replacing accordingly each factor in \eqref{FormulPathsCompleteTree} yields the desired formula \eqref{EnumComp}.
		\item   The proof of \eqref{EnumIncomp} follows from very similar arguments, using instead a binomial reproduction law  $\xi$ of parameter $\sigma$ and $\alpha$ (that is $\xi(k)= \alpha^{\sigma-k}\beta^{k}{\sigma \choose k}\mathbf 1_{k\leq \sigma}$) and displacement laws defined for every $j\leq \sigma$ by:
		\begin{align*}
			\eta^{(j)}=\frac{1}{|A^j_{\mathcal S}|}\sum_{\underline w\in A^j_{\mathcal S}} \delta_{\underline w},
		\end{align*}
		where $A_j^{\mathcal S}$ denotes the set of  all strictly increasing $j$-tuples of elements of $\mathcal S$. If $t$ is a canonically $\mathcal S$-embedded incomplete tree, we then have
		\begin{align*}
			\Pi_0(t)= \alpha^{(|t|+1)\sigma-|t|}\beta^{|t|}.
		\end{align*}
		It follows that
		\begin{align*}
		\Pi_0\left(e_{i, j}^T= x_{i, j}, \ \forall (i, j)\in \Delta_{\mathcal S}\right)=\alpha^{(N+1)\sigma-N}\beta^{N}\mathrm{Card}\left(\mathcal T^{\mathcal S}_{\text{incomp.}}(\underline x)\right).
		\end{align*}
        Moreover, in this case $\eta_*$ is just the law of $(Y_{-1}, Y_0, Y_1)$ where $Y_s=0$ if $s\notin \mathcal S$ and $(Y_s)_{s\in \mathcal S}$ are independent Bernoulli($\beta$) random variables. In particular, for every $n\geq 1$ and every $y_{-1}, y_0, y_1\leq n$, we have:
		\begin{align*}
			f^{\eta_*}_n(y_{-1}, y_0, y_1)= \1_{y_i=0 , \forall i\notin \mathcal S}  \prod_{i\in\mathcal S}\alpha^{n-y_i}\beta^{y_i}{n\choose y_i}.
		\end{align*}
		Then \eqref{EnumIncomp} follows from an application of \eqref{FormulPaths}, replacing each factor accordingly.
		\item To prove \eqref{EnumMotzkin}, we use the law $\xi=\frac{1}{3}(\delta_0+\delta_1+\delta_2)$ and the displacement laws $\eta^{(1)}=\delta_0$ and $\eta^{(2)}= \delta_{(-1, 1)}$ (cf \cite[Example 10.5]{Janson}). If $t$ is a canonically embedded Motzkin tree, we have:
		\begin{align*}
			\Pi_0(t)= 3^{-|t|-1}.
		\end{align*}
		This implies that
		\begin{align*}
				\Pi_0\left(e_{i, j}^T= x_{i, j}, \ \forall (i, j)\in \Delta\right)=3^{-N}\mathrm{Card}\left(\mathcal T_{\text{Motzkin}}(\underline x)\right)
		\end{align*}It is straightforward to check for this model that for every $n\geq 1$ and every $y_{-1}, y_0, y_1\leq n$ with $y_0+y_1\leq n$, we have:
		\begin{align*}
			f^{\eta_*}_n(y_{-1}, y_0, y_1)= \1_{y_{-1}=y_1}3^{-n}{n\choose y_{1}, y_0, n-y_0-y_1}.
		\end{align*}
		Then \eqref{EnumMotzkin} follows from an application of \eqref{FormulPaths}, replacing each factor accordingly.
	\end{enumerate}
\end{proof}

\section*{Acknowledgments}
The authors would like to thank  Igor Kortchemski, Jean-François Le Gall and Cyril Marzouk for their valuable comments and careful reading of earlier versions of this manuscript. This work was partially supported by the ERC Consolidator Grant 101087572, entitled “SuPerGRandMa”.

		\begin{small}
			\printbibliography
		\end{small}

\end{document}